# CONFIDENCE BANDS IN DENSITY ESTIMATION

By Evarist Giné and Richard Nickl

*University of Connecticut and University of Cambridge*

Given a sample from some unknown continuous density $f : \mathbb{R} \to \mathbb{R}$, we construct adaptive confidence bands that are honest for all densities in a "generic" subset of the union of $t$-Hölder balls, $0 < t \leq r$, where $r$ is a fixed but arbitrary integer. The exceptional ("nongeneric") set of densities for which our results do not hold is shown to be nowhere dense in the relevant Hölder-norm topologies. In the course of the proofs we also obtain limit theorems for maxima of linear wavelet and kernel density estimators, which are of independent interest.

**1. Introduction.** Let $X_1, \ldots, X_n$ be i.i.d. random variables with uniformly continuous density $f : \mathbb{R} \to \mathbb{R}$, and let $\hat{f}_n$ be some density estimator for $f$. A natural loss function to assess the statistical performance of $\hat{f}_n$ is sup-norm loss $d_\infty(\hat{f}_n, f) = \sup_x |\hat{f}_n(x) - f(x)|$: it gives a clear geometric interpretation of the estimation error, suggesting heuristically the existence of a "band" around $\hat{f}_n$ that shrinks at rate $d_\infty(\hat{f}_n, f)$ and contains $f$ with probability close to one.

Classical methods to construct confidence bands in density estimation—for example, the ones based on extreme value theory in Smirnov (1950) for histogram estimators and in Bickel and Rosenblatt (1973) for kernel estimators—require $f$ to satisfy stringent differentiability assumptions, and, more importantly, are based on *a priori knowledge* of the degree of smoothness of $f$. Recent developments in adaptive function estimation show that one can find purely data driven estimators $\hat{f}_n$ such that $d_\infty(\hat{f}_n, f)$ achieves the minimax-optimal rate of convergence $r_n(t) = (n/\log n)^{-t/(2t+1)}$ for estimating a density $f$ in a given $t$-Hölder ball [see Giné and Nickl (2009a, 2009b, 2010) in the i.i.d. density model on $\mathbb{R}$ and Goldenshluger and Lepski (2009) in the Gaussian white noise model]. The question then arises as to how one can take advantage of these adaptive rate of convergence results for statistical inference, in particular for the construction of "adaptive" confidence bands.

Let us phrase our problem more precisely: The results in Giné and Nickl (2009a, 2009b, 2010) show that the natural class of densities over which the minimax-optimal rate of convergence $r_n(t)$ can be achieved in $d_\infty$-loss in an adaptive way is

$$\mathcal{P}^{\text{all}} := \mathcal{P}^{\text{all}}(r, L) = \bigg\{ f : \mathbb{R} \to \mathbb{R} \text{ is a probability}$$
$$\text{density contained in } \bigcup_{0 < t \leq r} \Sigma(t, L) \bigg\},$$

where $\Sigma(t, L)$ is a ball of radius $L$ in the usual Hölder space on $\mathbb{R}$ and where the integer $r$ measures the degree of "regularity" of the kernel or wavelet basis used. Given $\alpha > 0$ and a family of densities $\mathcal{P} \subseteq \mathcal{P}^{\text{all}}$, an *honest confidence band* over the interval $[a, b]$ is a family of random intervals $C_n(y) := C_n(y, \alpha)$, $y \in [a, b]$, such that the asymptotic coverage inequality

$$(1.1) \qquad \liminf_n \inf_{f \in \mathcal{P}} \Pr_f(f(y) \in C_n(y) \text{ for all } y \in [a, b]) \geq 1 - \alpha$$

holds, and, following Cai and Low (2004), we shall say that $C_n(y)$ is *adaptive* if for every $t, \varepsilon > 0$ there exists $L'$ finite such that [$\ell(I)$ denoting the length of the interval $I$]

$$(1.2) \qquad \sup_{f \in \Sigma(t, L) \cap \mathcal{P}} \Pr_f \Big( \sup_y \ell(C_n(y)) \geq L' \tilde{r}_n(t) \Big) < \varepsilon,$$

where $\tilde{r}_n(t)$ equals $r_n(t)$, possibly inflated by a multiplicative logarithmic penalty.

It follows, on the one hand, from results in Low (1997) that confidence bands that are simultaneously adaptive and honest *do not exist* for $\mathcal{P} = \mathcal{P}^{\text{all}}$. *On the other hand we shall show that honest and adaptive confidence bands do exist over "generic" subsets $\mathcal{P}$ of $\mathcal{P}^{\text{all}}$.* The subset $\mathcal{P}$ of $\mathcal{P}^{\text{all}}$ for which our results hold is "generic" in the following sense:

(A) *$\mathcal{P}$ contains all "smooth" densities in $\mathcal{P}^{\text{all}}$, that is, all densities in $\mathcal{P}^{\text{all}}$ that are $r$-times differentiable on $\mathbb{R}$;*

(B) *The minimax rate of convergence over $\Sigma(t, L) \cap \mathcal{P}$ is the same as the one over $\Sigma(t, L) \cap \mathcal{P}^{\text{all}}$;*

(C) *The class of densities excised from $\mathcal{P}^{\text{all}}$ is "negligible" in the sense that the set $\mathcal{P}^{\text{all}} \setminus \mathcal{P}$ is "topologically small."*

Roughly speaking "topologically small" will mean that, given any $t > 0$, $\mathcal{P}$ can be chosen so large that the exceptional set $\mathcal{P}^{\text{all}} \setminus \mathcal{P}$ contains no (given) ball of $\Sigma(t,L) \cap \mathcal{P}^{\text{all}}$. If one relaxes the uniformity (or "honesty") requirement in (1.1) for the sake of illustration, our results will imply that $\mathcal{P}$ can be chosen so large that $\mathcal{P}^{\text{all}} \setminus \mathcal{P}$ is *nowhere dense* in the (relative) Hölder-norm topology in $\Sigma(t,L) \cap \mathcal{P}^{\text{all}}$ (again for every $t$). It should furthermore be noted that, although we state (B) separately, it typically follows from (C).

We will construct fully-data-driven adaptive (nonlinear) estimators $\hat{f}_n$ based on either wavelets or convolution kernels, and prove, uniformly over such a "generic" set $\mathcal{P}$, a "Smirnov–Bickel–Rosenblatt"-type limit theorem; for any bounded interval $[a,b]$,

$$(1.3) \qquad \hat{A}_n \left( \sup_{y \in [a,b]} \left| \frac{\hat{f}_n(y) - f(y)}{\hat{\sigma}_n \sqrt{\hat{f}_n(y)}} \right| - \hat{B}_n \right) \xrightarrow{d} Z$$

as $n \to \infty$ where $Z$ is a Gumbel random variable and where the random (but known) constants $\hat{\sigma}_n$, $\hat{A}_n$, $\hat{B}_n$ have the right stochastic order to obtain a confidence band from (1.3) that shrinks, up to a logarithmic penalty, at the minimax rate $r_n(t)$ of estimation. The estimator we propose is of "Lepski"-type and not difficult to implement. It is in principle possible to replace the interval $[a,b]$ by $\mathbb{R}$, by using suitable weight functions and techniques from Giné, Koltchinskii and Sakhanenko (2004), but this comes at the expense of much more technical proofs, so we abstain from it. See Section 3 for the exact statements of our results.

There has been substantial and deep recent work about the connection between confidence sets and rates of convergence of adaptive estimators. As mentioned above, Low (1997) shows some limitations for pointwise confidence intervals in density estimation. Our "generic" conditions circumvent his "pathologies." The paper closest to the present one is Picard and Tribouley (2000) where pointwise adaptive confidence intervals are constructed in regression and Gaussian white noise. Our proof strategy is partially inspired by theirs, and our Condition 3 is somewhat similar to their condition $H_s(M, x_0)$ which, however, is less "generic" in the pointwise setup. Cai and Low (2004) develop a general theory for pointwise confidence intervals, which can be conceptually (but not directly) related to the sup-norm case. Genovese and Wasserman (2008) revisit the negative results by Low (1997) in the framework of regression. They suggest that valid confidence sets are possible if the usual notion of "coverage" is replaced by "surrogate coverage," but it is not clear yet how "generic" this restriction is. There is also a remarkable literature on confidence sets in $L^2$-loss where the theory is somewhat different to the sup-norm/pointwise case, although the general message that "adaptive rates of convergence" do *not* simply translate into "adaptive confidence sets" is unchanged. We refer to Li (1989), Beran and Dümbgen (1998), Hoffman and Lepski (2002), Juditsky and Lambert-Lacroix (2003), Baraud (2004), Genovese and Wasserman (2005), Cai and Low (2006) and Robins and van der Vaart (2006). Another interesting approach is based on imposing qualitative shape constraints on the function to be estimated. Here some positive results are possible; we refer to Hengartner and Stark (1995), Dümbgen (2003) and Davies, Kovac and Meise (2009).

Most of the above literature is set in the Gaussian white noise model, but we prefer to derive our results in the i.i.d. density model, mostly for two reasons: First the asymptotic equivalence of white noise to density estimation only holds under quite restrictive assumptions on the underlying density; in particular we are interested in the low regularity case $t \leq 1/2$ as well. Second the problem of estimating a continuous density on $\mathbb{R}$ carries some specific structure that should be taken into account; such a density cannot be constant everywhere; neither can differentiable densities on $\mathbb{R}$ have derivatives that are everywhere zero, facts that play a role in the verification of some of our conditions.

The limit (1.3) is based on conditions that require certain centered linear wavelet or kernel estimators to satisfy a "Smirnov–Bickel–Rosenblatt"-type limit theorem, uniformly in the underlying density $f$. While these results can be obtained, as we show, for convolution kernel estimators along the lines of Bickel and Rosenblatt (1973), using refinements from Giné, Koltchinskii and Sakhanenko (2004), results of this type *do not exist at the moment* for wavelet density estimators. It turns out that a reduction to Gaussian processes similar to the one for kernel estimators first proposed by Bickel and Rosenblatt (1973) can be proved for wavelets as well (see Proposition 5), but the resulting Gaussian process, which equals

$$Y(t) = \int_{\mathbb{R}} K(x,y) \, dW(y),$$

where $K$ is the wavelet projection kernel and $W$ is Brownian motion, turns out to be *nonstationary*, so that the classical extreme value theory for stationary Gaussian processes [Leadbetter, Lindgren and Rootzén (1983)] *does not apply here*. However, these Gaussian processes are *cyclostationary* [meaning that the covariance function $r(t, t+v)$ is periodic in $t$ with the same period for all $v$], and the extreme value theory for these and related processes has recently attracted some interest in the literature [see Konstant and Piterbarg (1993), Piterbarg and Seleznjev (1994), Hüsler (1999) and Hüsler, Piterbarg and Seleznjev (2003)]. Using these techniques and wavelet theory we can establish, as a first step, limit theorems for suprema of centered wavelet density estimators based on *Battle–Lemarié* wavelets [which can be computed numerically as spline projections, cf. Giné and Nickl (2010)]. We believe that this proof strategy should also work for other wavelet bases, but we currently do not have enough knowledge about the analytical properties of the covariance functions of the processes $Y(t)$ in the case of, for instance, Daubechies wavelets, to succeed in doing so. This remains an open problem.

We finally remark that the results in this article are clearly of an asymptotic (and hence "theoretical" nature); they show that adaptive confidence bands are possible for large sample sizes and in a certain "generic" sense, but we do not advocate the use of our bands in practice without a thorough investigation of their finite sample properties.

## 2. Basic notation and definitions.

2.1. *Wavelets and function spaces.* For a function $H:M \to \mathbb{R}$ we shall denote by $\|H\|_M$ the quantity $\sup_{m \in M} |H(m)|$, but we shall write $\|H\|_\infty := \|H\|_\mathbb{R}$. Denote further by $C(\mathbb{R})$ the space of bounded continuous functions on $\mathbb{R}$ normed by $\|\cdot\|_\infty$.

We next define the function spaces that will be at the heart of our statistical problem. It is convenient to define them in terms of wavelets; more classical equivalent definitions can be found in the literature (see Remark 1 below). Throughout this paper we shall use the by now standard wavelet theory, we refer, for example, to the monograph Härdle et al. (1998) in what follows—for an excellent treatment of the statistically most relevant materials. In particular we shall say the scaling function $\phi$ of a multiresolution analysis is *s-regular* if $\phi$ is $s$-times weakly differentiable, and, for $0 \le \alpha \le s$, $D^\alpha \phi$ satisfies $|D^\alpha \phi(x)| \le c_1 \lambda^{c_2|x|}$ for (almost) every $x \in \mathbb{R}$, some $c_1, c_2 > 0$ and some $0 < \lambda < 1$.

DEFINITION 1. Let $0 < t < s$, $t \in \mathbb{R}$, $s \in \mathbb{N}$. Let $\phi$ be a scaling function that is $s$-regular; let $\psi$ be the associated mother wavelet and denote by $\alpha_k(f)$ and $\beta_{lk}(f)$, $k \in \mathbb{Z}$, $l \in \mathbb{N}$, the associated wavelet coefficients of the function $f:\mathbb{R} \to \mathbb{R}$. The Hölder–Zygmund space $\mathcal{C}^t(\mathbb{R})$ is defined as the set of functions

$$\mathcal{C}^t(\mathbb{R}) := \Big\{ f \in C(\mathbb{R}) : \|f\|_{t,\infty} := \sup_{k \in \mathbb{Z}} |\alpha_k(f)|$$
$$+ \sup_{l \ge 0} \sup_{k \in \mathbb{Z}} |2^{l(t+1/2)} \beta_{lk}(f)| < \infty \Big\}.$$

We should note that this definition is independent of the wavelet basis used: any wavelet basis of regularity $s > t$ generates the same space.

REMARK 1. It is a standard result in wavelet theory [e.g., Chapter 6.4 in Meyer (1992)] that $\mathcal{C}^t(\mathbb{R})$ is equal, with equivalent norms, to the classical Hölder–Zygmund spaces, defined as follows: For $0 < t < 1$ define $C^t(\mathbb{R})$ to be the space of functions $f \in C(\mathbb{R})$ for which $\|f\|'_{t,\infty} := \|f\|_\infty + \sup_{x \ne y, x, y \in \mathbb{R}} (|f(x) - f(y)|/|x - y|^t)$ is finite. For noninteger $t > 1$ the space $C^t(\mathbb{R})$ is defined by requiring $D^{[t]} f$ of $f \in C(\mathbb{R})$ to exist and to be contained in $C^{t-[t]}(\mathbb{R})$. The Zygmund class $C^1(\mathbb{R})$ is defined by requiring $|f(x+y) + f(x-y) - 2f(x)| \le C|y|$ for all $x, y \in \mathbb{R}$, some $0 < C < \infty$ and $f \in C(\mathbb{R})$, and the case $m < t \le m+1$ follows by requiring the same condition on the $m$th derivative of $f$. It is then also clear that $C^m(\mathbb{R})$ contains the spaces of $m$-times continuously differentiable functions with $m$ bounded derivatives. We remark finally that $C^t(\mathbb{R})$ is a special case of the scale of Besov spaces, namely $B^t_{\infty\infty}(\mathbb{R})$.

2.2. *Density estimation using convolution kernels or wavelets.* Let $X_1, \ldots, X_n$ be i.i.d. random variables with common law $P$ and density $f$ on $\mathbb{R}$, and split the sample into two parts, $S_1$ and $S_2$, each of size $n_1$ and $n_2$, respectively, in such a way that $n_1/n_2$ is bounded away from zero and infinity as $n \to \infty$. Denote by

$$P_{n_1} = \frac{1}{n_1} \sum_{i=1}^{n_1} \delta_{X_i} \quad \text{and} \quad P_{n_2} = \frac{1}{n_2} \sum_{i=1}^{n_2} \delta_{X_{n_1+i}}$$

the empirical measures associated with the first and the second subsample. We take the $\{X_i\}$'s to be coordinate projections of the infinite product probability space $\mathbb{R}^\mathbb{N}$ with its product sigma-algebra and denote by $\Pr_f$ the product probability measure on this space.

We will consider two types of preliminary linear estimators: Define first the classical kernel density estimator based on the sample $S_v$, $v = 1, 2$, namely

$$\frac{1}{h} \int_\mathbb{R} K\left(\frac{y-x}{h}\right) dP_{n_v}(x), \qquad y \in \mathbb{R},$$

where $K : \mathbb{R} \to \mathbb{R}$ is a kernel and $h > 0$ is some bandwidth. An alternative estimator is based on a wavelet projection: If $\phi$ is a scaling function (father wavelet) and $\psi$ the associated (mother) wavelet, then the linear wavelet estimator based on the sample $S_v$, $v = 1, 2$ is

$$2^j \int_\mathbb{R} K(2^j y, 2^j x) \, dP_{n_v}(x)$$

$$= \sum_k \hat{\alpha}_k(v) \phi(y-k) + \sum_{l=0}^{j-1} \sum_k \hat{\beta}_{lk}(v) \psi_{lk}(y), \qquad y \in \mathbb{R}, j \in \mathbb{N},$$

where $K(y, x) = \sum \phi(y-k)\phi(x-k)$ is the wavelet projection kernel, $\psi_{lk}(x) = 2^{l/2} \psi(2^j x - k)$ and where the empirical wavelet coefficients are

$$\hat{\alpha}_k(v) = \int_\mathbb{R} \phi(x-k) \, dP_{n_v}(x), \qquad \hat{\beta}_{lk}(v) = \int_\mathbb{R} \psi_{lk}(x) \, dP_{n_v}(x).$$

To unify the notation for both estimators we convert the bandwidth $h$ into $2^{-j}$ so that the *kernel-type* density estimator is given by

(2.1) $$f_{n_v}(y, j) = 2^j \int_\mathbb{R} K(2^j y, 2^j x) \, dP_{n_v}(x),$$

where $K(y,x)$ is either the wavelet projection kernel or the convolution kernel $K(y-x)$. In this way the estimator is defined also for noninteger $j$. We will use the convention $K_j(y,x) = 2^j K(2^j y, 2^j x)$, and we denote the expectation of $f_{n_v}(y,j)$ by

$$(2.2) \qquad E f_{n_v}(y,j) = \int_{\mathbb{R}} K_j(y,x) f(x)\, dx = K_j(f)(y).$$

We shall make the following standard assumption on the kernel $K$ which will have to be strengthened for some results.

CONDITION 1. Let $r \in \mathbb{N}$, $r \geq 1$. Suppose one of the following conditions is satisfied:

(a) (convolution kernel) let $K(x,y) = K(x-y)$ where $K: \mathbb{R} \to \mathbb{R}$ is symmetric, integrable, of bounded variation and integrates to one. Assume furthermore that $\int_{\mathbb{R}} K(u) u^l = 0$ for $l = 1, \ldots, r-1$ (vacuous if $r=1$) as well as $\int_{\mathbb{R}} |K(u)| |u|^r\, du < \infty$;

(b) (wavelet kernel) (i) let $K(x,y) = \sum_k \phi(x-k)\phi(y-k)$ where $\phi$ is a scaling function that is of bounded variation, compactly supported and either $\phi$ is $(r-1)$-regular or $\psi$ satisfies $\int_{\mathbb{R}} \psi(u) u^l = 0$ for every $0 \leq l \leq r-1$;

(c) (wavelet kernel) (ii) let $K(x,y) = \sum_k \phi_r(x-k)\phi_r(y-k)$ where $\phi_r$ is the *Battle–Lemarié* scaling function (defined in Section 4.3.1).

## 3. Adaptive confidence bands.

3.1. *Estimate of the resolution level.* We use the sample $S_2$ to choose the resolution level $j$. For integers $r \geq 1$, $n_2 > 1$, choose integers $j_{\min} := j_{\min,n}$ and $j_{\max} := j_{\max,n}$, $0 < j_{\min} < j_{\max}$, such that

$$2^{j_{\min}} \simeq \left(\frac{n_2}{\log n_2}\right)^{1/(2r+1)}, \qquad 2^{j_{\max}} \simeq \left(\frac{n_2}{(\log n_2)^4}\right),$$

$$\mathcal{J} := \mathcal{J}_n = [j_{\min}, j_{\max}] \cap \mathbb{N}.$$

We note in advance that $j_{\min}$ is the resolution level we would choose if we knew that the unknown density $f$ is $r$-times continuously differentiable, and we are not trying to adapt to densities smoother than this. (This does not mean that we rule out densities that are very smooth; it just means that we live with a "nonadaptive" rate of convergence in these cases.) On the other hand, $j_{\max}$ is the resolution level that just produces uniform consistency of the linear estimator $f_{n_2}(j_{\max})$ if $f$ is bounded and uniformly continuous. So our problem is to adapt to the unknown smoothness $t$ of $f$ where $t$ varies between 0 and $r$.

Our data-driven choice $\hat{j}_n$ for the resolution level is of "Lepski-type" and is based on the subsample $S_2$; namely

$$(3.1) \quad \hat{j}_n = \min\left\{j \in \mathcal{J} : \|f_{n_2}(j) - f_{n_2}(l)\|_\infty \leq M\sqrt{\frac{2^l l}{n_2}} \ \forall l > j, l \in \mathcal{J}\right\},$$

where $M = M'\sqrt{\|f\|_\infty \vee 1}$ with $M' := M'(K)$, a constant that depends only on $K$. We discuss in Remark 2 below the choice of $M'$ as well as how one can circumvent having to know $\|f\|_\infty$ in practice.

A remark on the choice of $\hat{j}_n$ is in order. Our proofs will imply, for $0 < t < r$, the adaptive global (minimax-optimal) risk bound

$$\sup_{f:\|f\|_{t,\infty} \leq L} E \sup_{y \in \mathbb{R}} |f_{n_1}(y, \hat{j}_n) - f(y)| = O\left(\left(\frac{\log n}{n}\right)^{t/(2t+1)}\right).$$

To make inferential use of this result we will use the estimator $f_{n_1}(\hat{j}_n)$ (in fact a slight modification of it) as the center of a confidence band for the unknown density $f$ on the interval $[a, b]$. Under certain assumptions on $f$ our confidence band will be shown to be both honest and adaptive for *arbitrary* bounded intervals $[a, b]$ (although we prove our result, w.l.o.g., only for $[a, b] = [0, 1]$). If one starts with a fixed interval $[a, b]$, one may alternatively try to choose the resolution level $\hat{j}_n$ above depending only on values of $f_{n_2}(j), f_{n_2}(l)$ on $[a, b]$. This is not our approach here, however; we want to construct a single estimator that is *globally* adaptive, and find honest confidence bands on arbitrary intervals $[a, b]$ for it. The important question of spatial adaptation is not addressed in the present paper.

3.2. *The main assumptions.* For the main theorem below, we will need some conditions that we state now. The first condition is stochastic in nature and is about an exact limit theorem for the maximum deviations of the centered linear estimator. Define

$$(3.2) \qquad c(K) = \sqrt{\sup_x \int_{\mathbb{R}} K^2(x, y)\, dy},$$

and let $A(l) := A(l, K), B(l) := B(l, K)$ be real-valued functions defined on $\mathbb{N}$, depending only on $K$ and such that $A(l) \simeq l^{1/2} \simeq B(l)$. For $n, l \in \mathbb{N}, x \in \mathbb{R}$, define

$$T_n(l, x, f, K)$$
$$:= \left|\Pr_f\left\{\sqrt{\frac{n_1}{2^l}} \sup_{y \in [0,1]} \left|\frac{f_{n_1}(y, l) - Ef_{n_1}(y, l)}{c(K)\sqrt{f(y)}}\right| \leq \frac{x}{A(l)} + B(l)\right\} - e^{-e^{-x}}\right|.$$

For $F = [F_1, F_2]$ with $F_1 < 0 < 1 < F_2$ and $\delta, \alpha > 0$ define further the class of densities

$$
\begin{aligned}
\mathcal{D} &= \mathcal{D}(\alpha, D, \delta, F) \\
&= \left\{ f : \mathbb{R} \to \mathbb{R}, \int_{\mathbb{R}} f = 1, f \geq 0 \text{ on } \mathbb{R}, f \geq \delta \text{ on } F, \|f\|_{\alpha, \infty} \leq D \right\}.
\end{aligned}
\tag{3.3}
$$

To avoid triviality we shall only consider combinations of $\alpha, D, \delta, F$ such that $\mathcal{D}$ is nonempty, and given $\alpha, D$ we shall say that $\delta, F$ are "admissible" if $\mathcal{D}$ is nonempty.

CONDITION 2. Assume that for every $x \in \mathbb{R}$, every $0 < \alpha < \infty$, $0 < D < \infty$, $l_n \in \mathcal{J}_n$ and every admissible $\delta > 0$, $F$ we have, as $n \to \infty$, that

$$\sup_{f \in \mathcal{D}(\alpha, D, \delta, F)} |T_n(l_n, x, f, K)| \to 0.$$

Verifying this condition is a nontrivial problem in itself, and we discuss this in detail in Section 3.4. We also need the following condition on the underlying density.

CONDITION 3. Suppose $f \in \mathcal{C}^t(\mathbb{R})$ for some $t > 0$ and that there exist positive finite constants $b_1 \leq b_2$ and a positive integer $j_0$ such that for every integer $j \geq j_0$,

$$b_1 2^{-jt} \leq \|K_j(f) - f\|_\infty \leq b_2 2^{-jt}.$$

Note that the upper bound is standard and can be shown to follow from $f \in \mathcal{C}^t(\mathbb{R})$ for $t < r$ and $K$ satisfying Condition 1 [cf. Theorem 9.3 in Härdle et al. (1998)]. For the uniformity results below it is convenient to require the upper bound also in the boundary case $t = r$ [which does not necessarily follow from just $f \in \mathcal{C}^r(\mathbb{R})$, but holds, for example, for $r$-times differentiable functions with $r$ bounded derivatives; cf. Theorem 8.1 in Härdle et al. (1998)]. The lower bound on the error of approximation of $f$ by $K_j(f)$ is a crucial assumption, and we refer to Section 3.5 for a detailed discussion of this condition.

In the construction of confidence intervals, it is well known [since Bickel and Rosenblatt (1973); see also Hall (1992) and Picard and Tribouley (2000)] that one should "undersmooth." In the classical case of convolution kernels this means that we should decrease the bandwidth $2^{-\hat{j}_n}$ to $2^{-\hat{j}_n - u_n}$ and use the function $f_{n_1}(y, \hat{j}_n + u_n)$ as the center of the band where $u_n$ is some sequence of positive numbers. In the context of wavelets this means that we should add a block of empirical wavelet coefficients at resolutions $\hat{\hat{j}}_n \leq j < \hat{\hat{j}}_n + u_n$ to our estimator so that the function

$$(3.4) \qquad f_{n_1}(y, \hat{\hat{j}}_n + u_n) = f_{n_1}(y, \hat{\hat{j}}_n) + \sum_{l=\hat{\hat{j}}_n}^{\hat{\hat{j}}_n + u_n - 1} \sum_k \hat{\beta}_{lk}(1) \psi_{lk}(y)$$

is the center of the confidence band.

CONDITION 4. Let $u_n$ be a sequence of positive integers such that $2^{u_n} \simeq (\log n)^2$.

3.3. *The main result.* Let $\hat{\hat{j}}_n$ be the data-driven resolution level from (3.1); recall the constants from Condition 2, and define

$$\hat{\sigma}_n = \sqrt{\frac{2^{\hat{\hat{j}}_n + u_n}}{n_1}}, \qquad \hat{A}_n = A(\hat{\hat{j}}_n + u_n), \qquad \hat{B}_n = B(\hat{\hat{j}}_n + u_n).$$

If $c(K)$ is as in (3.2), then the size of the band around $f_{n_1}(y, \hat{\hat{j}}_n + u_n)$ will be twice

$$(3.5) \qquad s_n(y, x) = \hat{\sigma}_n c(K) \sqrt{f_{n_1}(y, \hat{\hat{j}}_n + u_n)} \left( \frac{x}{\hat{A}_n} + \hat{B}_n \right),$$

and we note that this quantity can be shown to be eventually positive for *every* $x \in \mathbb{R}$, but for fixed $n$ we implicitly assume that $x$ is large enough so that $s_n(y, x) > 0$. We emphasize that $s_n(x)$ is completely data-driven [except for the dependence on $\|f\|_\infty$ through (3.1) discussed in Remark 2]. The confidence band we propose for $f$ is

$$(3.6) \qquad \begin{aligned} C_n(x, y) := [&f_{n_1}(y, \hat{\hat{j}}_n + u_n) - s_n(y, x), \\ &f_{n_1}(y, \hat{\hat{j}}_n + u_n) + s_n(y, x)], \qquad x, y \in \mathbb{R}, \end{aligned}$$

and the probability of inferential interest is, for $x \in \mathbb{R}$,

$$\Pr_f \{ f(y) \in C_n(x, y) \text{ for every } y \in [0, 1] \}.$$

As mentioned before, we restrict ourselves here to the interval $[0, 1]$, but any bounded set $[a, b]$ is possible as long as $f$ is bounded away from zero on a neighborhood of $[a, b]$.

Our adaptation result will be shown to hold for densities satisfying Condition 3 with $t \in (0, r]$, where $r$ is the regularity of the kernel $K$ from Condition 1. Moreover, the result will be uniform for $t$ in any compact subset of $(0, r]$.

To describe exactly the set of densities over which our results hold uniformly, define first, for fixed $0 < \eta \leq r$, $0 < b < \infty$, $0 < b_1 \leq b_2 < \infty$ and $j_0 \in \mathbb{N}$,

$$\mathcal{P}(\eta, r, b, b_1, b_2, j_0)$$
$$:= \bigcup_{\eta \leq t \leq r} \{f \in \mathcal{C}^t(\mathbb{R}) : \|f\|_{\eta,\infty} \leq b, b_1 2^{-jt} \leq \|K_j(f) - f\|_\infty$$
$$\leq b_2 2^{-jt} \ \forall j \geq j_0\}.$$

This class is just the union over $t \in [\eta, r]$ of functions that satisfy Condition 3 for the given $t$ (and constants $b_1, b_2, j_0$), and that is also contained in a fixed ball of $\mathcal{C}^\eta(\mathbb{R})$. We assume implicitly that $b$ is large enough so that this class is nonempty.

Let then $\mathcal{D}(\eta, b, \delta, F)$ be the set from (3.3) where $\delta, F$ are admissible, and define

(3.7) $\quad \mathcal{P} := \mathcal{P}(\eta, r, b, b_1, b_2, j_0, \delta, F) = \mathcal{P}(\eta, r, b, b_1, b_2, j_0) \cap \mathcal{D}(\eta, b, \delta, F).$

This set simply consists of densities that are in $\mathcal{P}(\eta, r, b, b_1, b_2, j_0)$ and that are also bounded away from zero on $F$. It is easy to see that this set is nonempty, and we shall discuss this class of densities in detail in Section 3.5. Clearly, for every $f \in \mathcal{P}$ there exists a unique $t := t(f)$ for which Condition 3 is satisfied.

THEOREM 1. *Let $f_{n_1}(y, l)$ be the estimator from (2.1) with $K$ satisfying Condition 1 for some $r \geq 1$. Let further $\hat{j}_n$ be defined as in (3.1), and let $c(K)$ be as in (3.2). Assume that Conditions 2 and 4 are satisfied for $f_{n_1}(y, l)$ and $u_n$, respectively. Then we have for every $x \in \mathbb{R}$, $b > 0$, $0 < b_1 \leq b_2 < \infty, j_0 \in \mathbb{N}$, $0 < \eta \leq r$ and every admissible $\delta > 0, F = [F_1, F_2]$ satisfying $F_1 < 0 < 1 < F_2$ that*

$$\sup_{f \in \mathcal{P}(\eta, r, b, b_1, b_2, j_0, \delta, F)} \left| \Pr_f \left\{ \hat{A}_n \left( \sup_{y \in [0,1]} \left| \frac{f_{n_1}(y, \hat{j}_n + u_n) - f(y)}{c(K)\hat{\sigma}_n \sqrt{f_{n_1}(y, \hat{j}_n + u_n)}} \right| - \hat{B}_n \right) \leq x \right\} \right.$$
$$\left. - e^{-e^{-x}} \right|$$

*converges to 0 as $n \to \infty$. Furthermore, for every $\varepsilon > 0$ there exists a constant $L$ such that, for every $n \in \mathbb{N}$, ($t = t(f)$)*

(3.8) $\quad \sup_{f \in \mathcal{P}(\eta, r, b, b_1, b_2, j_0, \delta, F)} \Pr_f \{\hat{\sigma}_n \geq Ln^{-t/(2t+1)}(\log n)^{-1/2(2t+1)} 2^{u_n/2}\} \leq \varepsilon.$

In the proof, which is given in Section 4.4.1, we will show that

$$0 < \inf_{y \in [0,1]} f_n^{1/2}(y, \hat{j}_n + u_n) \leq \sup_{y \in [0,1]} f_n^{1/2}(y, \hat{j}_n + u_n) < L'$$

for some constant $L'$ on sets of probability approaching one, and the fraction in the above theorem has to be understood accordingly. Moreover, $A(l) \simeq \sqrt{l}$ (Condition 2), Condition 4 and $\hat{j}_n \in \mathcal{J}$ imply

$$\hat{A}_n = A_n(\hat{j}_n + u_n) \simeq \sqrt{\log n},$$

and likewise for $\hat{B}_n$. Combining this with Condition 4 we have the following:

COROLLARY 1. *Let the assumptions of Theorem 1 be satisfied; let $C_n(y,x)$ be the confidence band from (3.6), and let $\mathcal{P}$ be as in (3.7). Then, for every $x \in \mathbb{R}$,*

$$\sup_{f \in \mathcal{P}} |\Pr_f\{f(y) \in C_n(y,x) \ \forall y \in [0,1]\} - e^{-e^{-x}}|$$

*converges to zero as $n \to \infty$. Furthermore, this confidence band is adaptive: If $2s_n(y,x)$ is the length of $C_n(y,x)$ at $y \in [0,1]$, then for every $\varepsilon > 0$ there exists a constant $L$ such that, for every $n \in \mathbb{N}$ and every $x \in \mathbb{R}$, ($t = t(f)$)*

$$(3.9) \qquad \sup_{f \in \mathcal{P}} \Pr_f \left\{ \sup_{y \in [0,1]} s_n(y,x) \geq L(n/\log n)^{-t/(2t+1)} 2^{u_n/2} \right\} \leq \varepsilon.$$

REMARK 2. The definition of $\hat{j}_n$ in (3.1) involves two "unknown" quantities. The first is $\|f\|_\infty$, and in practice this can be replaced, for instance, by the estimate $\|f_{n_2}(j_{\max})\|_\infty$. All proofs go through for this data-driven choice as well [arguing as in the proof of Theorem 2 of Giné and Nickl (2009a)], but we abstain from this to reduce technicalities. Another question is how to select the constant $M'$. A concrete choice can be obtained from tracking the constants in the proof of Lemma 2. To obtain good constants one may use Rademacher-symmetrization in a similar vein as in Giné and Nickl (2010).

3.4. *Condition 2 and the asymptotic distributions of suprema of linear density estimators.* Since the results in this section only involve the sample $S_1$, we set $n_1 = n$. The prototypical result required in Condition 2 [without uniformity in $\mathcal{D}(\alpha, D, \delta, F)$] is for convolution kernel density estimators, and due to Bickel and Rosenblatt (1973). Their conditions are too stringent for our "adaptive" framework, but some refined methods from Giné, Koltchinskii and Sakhanenko (2004) can be used to verify Condition 2.

PROPOSITION 1. *If $K : \mathbb{R} \to \mathbb{R}$ satisfies Condition 1(a), is supported in $[-1,1]$ and is twice continuously differentiable on $\mathbb{R}$, then the kernel density estimator $f_n(y,l)$ from (2.1) satisfies Condition 2 with $A(l) = \sqrt{2(\log 2)l}$, $B(l)$ defined in (4.22), and $c(K) = \|K\|_2$.*

PROOF. Use Proposition 7 below. □

One is next led to ask whether an analogue of the classical Bickel–Rosenblatt theorem can be proved for the wavelet case as well. This problem has no simple solution in general; the easiest case being that of Haar wavelets which was already considered in Smirnov (1950). Let $f_n(y,l)$ be as in (2.1) where $\phi = 1_{[0,1)}$ which satisfies Condition 1(b) with $r = 1$.

PROPOSITION 2. *The Haar-wavelet density estimator satisfies Condition 2 with $A(l)$ and $B(l)$ defined in (4.21) and with $c(K) = 1$.*

PROOF. Use Proposition 6 below. □

The Haar-wavelet allows one to adapt only up to smoothness one, so it is of interest to verify Condition 2 for other wavelets that satisfy Condition 1 for $r \geq 2$. On the positive side we prove in Section 4.1 a Gaussian reduction argument for general wavelet estimators similar to the convolution kernel case. The resulting Gaussian processes are given by the stochastic integrals

$$Y(t) = \int_{\mathbb{R}} K(t,u)\, dW(u),$$

where $W$ is Brownian motion and where $K(x,y) = \sum_k \phi(x-k)\phi(y-k)$ is the wavelet projection kernel. On the negative side, these processes are *nonstationary* and therefore we cannot use the classical extreme value theory for stationary Gaussian processes [Leadbetter et al. (1983)] as Bickel and Rosenblatt (1973) and Giné, Koltchinskii and Sakhanenko (2004) did.

The theory for nonstationary Gaussian processes is more involved (see Section 4.2.3). In Section 4.3 we will prove that the wavelet density estimators based on Battle–Lemarié wavelets of degree $r \leq 4$ (which satisfy condition 1 for this $r$) do satisfy Condition 2. We believe that the condition holds for Battle–Lemarié wavelets of any degree, but our proof depends on specific computations that increase in complexity with the degree, and which we complete only for $r \leq 4$. See Remark 3 after the proof of Proposition 9 below for more discussion. The case $r = 1$ (Haar wavelet) is not repeated.

PROPOSITION 3. *The wavelet density estimator $f_n(y,l)$ from (2.1) based on Battle–Lemarié wavelets $\phi_r$ with $r \in [2,4]$ satisfies Condition 2 with $A(l)$, $B(l)$ and $c(K) = \sigma_r$ as in Propositions 8 and 9.*

3.5. *Condition 3 and the class $\mathcal{P}$.* As was mentioned in the introduction, the natural classes for adaptive density estimation in sup-norm loss are balls $\Sigma(t, L)$ of radius $L$ in $\mathcal{C}^t(\mathbb{R})$, where $0 < t < r$ (including the case $t = r$ if the upper bound in Condition 3 holds with $t = r$). Theorem 1 does not hold for $\bigcup_{0<t\leq r} \Sigma(t,L)$, but only for $\mathcal{P}$, and we want to discuss in detail this class in order to understand the restrictions imposed, mostly Condition 3. On the one hand, we recall from the introduction that a honest adaptive confidence band cannot exist for the full class $\bigcup_{0<t\leq r} \Sigma(t,L)$. On the other hand we shall show below that (i) any $r$-times differentiable density with bounded continuous derivatives is contained in $\mathcal{P}$ (for some $b_1, b_2$), so our confidence band is valid, and shrinks at rate $n^{-r/(2r+1)}$ (up to a logarithmic term) for very smooth densities; (ii) the minimax rate of convergence over $\mathcal{P} \cap \Sigma(t,L)$ is the same as the one over $\Sigma(t,L)$, and (iii) one cannot "generically" improve upon the class $\mathcal{P}$ in Theorem 1, at least in the following sense: *The set of densities that are contained in $\Sigma(t,L)$ but not in $\mathcal{P}$ contains no given ball of $\Sigma(t,L)$.* Exact statements require some more careful discussion.

Let us first remark that the mild requirement that the density $f$ is bounded away from zero on an interval whose interior contains $[0,1]$ is helpful in the verification of Condition 2. This condition could be avoided by using techniques from Giné, Koltchinskii and Sakhanenko (2004) but at the expense of considerably more technical proofs (that also lead to modified results).

The crucial restriction that we impose is Condition 3. Verification of the upper bound in that condition is standard and was already discussed immediately after Condition 3. The delicate part of the condition is the lower bound. We start with an informal discussion in the case where $K$ is a convolution kernel, and taking $K = 1_{[-1/2, 1/2]}$ for simplicity (so that $r=2$). In this case the quantity in this condition reduces to

$$(3.10) \qquad \|K_j(f) - f\|_\infty = \sup_{x \in \mathbb{R}} \left| \int_{-1/2}^{1/2} (f(x - u2^{-j}) - f(x))\, du \right|.$$

First, if the density $f$ is infinitely differentiable with bounded continuous derivatives, then this quantity is of order

$$2^{-2j-1} \sup_{x \in \mathbb{R}} \left| \int_{-1/2}^{1/2} u^2 D^2 f(x)\, du \right| + o(2^{-2j}) \geq b_1 2^{-2j}$$

for some $b_1 > 0$ and $j$ large enough since no such density on $\mathbb{R}$ can have a second derivative that is everywhere zero. The constant $b_1$ is even bounded away from zero uniformly in the set of all twice differentiable densities that are supported in a fixed compact interval $[a,b]$ [as is easily seen by expanding $0 = f(a)$ up to second-order around a point of maximum $x_0$, using also $f(x_0) \geq (b-a)^{-1}$]. If $K$ is a kernel of order $r$ but not of order $r+1$, then the same lower bound holds with $2^{-2j}$ replaced by $2^{-rj}$, and *we see that Condition 3 is then always satisfied with $t = r$ for very smooth densities $f$.* Similar remarks apply to wavelets.

Hence we have to consider the case where $f$ is *not* very smooth, say $f \in \mathcal{C}^t(\mathbb{R})$, but *not in* $\mathcal{C}^{t+\gamma}(\mathbb{R})$ for any $\gamma > 0$. For a given function $f$, one can call, in slight abuse of terminology,

$$t(f) := \sup\{t : f \in \mathcal{C}^t(\mathbb{R})\}$$

its "Hölder exponent." This exponent is generally not attained for a given function $f$, but before we address this issue let us continue with some special cases. Suppose for instance $f$ is infinitely differentiable except at $x_0$ where $f$ behaves locally as $|x - x_0|$ so that $f \in \mathcal{C}^1(\mathbb{R})$ but $f \notin \mathcal{C}^{1+\gamma}(\mathbb{R})$ for any $\gamma > 0$. This means that we would like to verify Condition 3 with $t = 1$. Indeed the integrand in (3.10), for $x = x_0$, equals $2^{-j}|u|$ so that again $\|K_j(f) - f\|_\infty \geq b_1 2^{-j}$. More generally we can rewrite the quantity in (3.10) as

$$2^{-jt} \sup_{x \in \mathbb{R}} \left| \int_{-1/2}^{1/2} |u|^t \frac{f(x - u2^{-j}) - f(x)}{(|u|2^{-j})^t} \, du \right|.$$

Now intuitively we would expect that $f$ in $\mathcal{C}^t(\mathbb{R})$ but *not in* $\mathcal{C}^{t+\gamma}(\mathbb{R})$ for any $\gamma > 0$ precisely means that $f$ attains the Hölder exponent $t$ and that, for some $x_0 \in \mathbb{R}$,

$$\frac{f(x_0 - v) - f(x_0)}{|v|^t}$$

is bounded away from zero or even has a nonzero limit as $v \to 0$ [where $|v|^t$ has to be replaced by $\text{sign}(v)|v|^t$ in the denominator]. Unfortunately this reasoning is too naive, and it is not difficult to see that $\mathcal{C}^t(\mathbb{R})$ contains functions that do not attain their Hölder exponent; in particular there exist functions in $\mathcal{C}^t(\mathbb{R}) \setminus (\bigcup_{s > t} \mathcal{C}^s(\mathbb{R}))$ for which the Hölder exponent is not $t$. However, one can show that such a pathology cannot occur for "quasi-every" function in $\mathcal{C}^t(\mathbb{R})$.

To be precise, let us recall that a property holds for "quasi-every" element in a metric space if the set of elements in this space that do not satisfy this property is *nowhere dense* (so in particular "meagre" in the sense of Baire categories). Recall further that a subset $F$ of a metric space is nowhere dense if the interior of its closure is empty, so in particular $F$ contains no open subset. For example, a classical result of Banach (1931) is that "quasi-every" function in the space of continuous functions on $[0, 1]$ (equipped with the sup-norm) is nowhere differentiable. This is sensible: for instance, any bounded set of equicontinuous functions is norm-compact in this space, and compact sets in infinite-dimensional normed linear spaces are always meagre.

Inspired by these ideas and recent results of Jaffard (1997, 2000), we can now state the main result of this subsection which says that the set of functions in the Banach space $\mathcal{C}^t(\mathbb{R})$ that do not satisfy the lower bound in Condition 3 is nowhere dense in this space.

PROPOSITION 4. *For $f \in \mathcal{C}^t(\mathbb{R})$, let $K_j(f)$ be as in (2.2) where $K(x,y) = \sum_k \phi(x-k)\phi(y-k)$ and where $\phi$ is $(r-1)$-regular, $r-1 > t$. Then the set of functions*

$$\mathcal{N}_t = \{f \in \mathcal{C}^t(\mathbb{R}) \colon \text{there do not exist } b_1 > 0$$
$$\text{and } j_0 \in \mathbb{N} \text{ s.t. } \|K_j(f) - f\|_\infty \geq b_1 2^{-jt} \ \forall j \geq j_0\}$$

*is nowhere dense in (the norm-topology of) the Banach space $\mathcal{C}^t(\mathbb{R})$.*

The proof can be found in Section 4.4.2. Note further that the proofs for convolution kernels and also for $r-1 \leq t \leq r$ are similar but more technical.

The question now arises as to how exactly this result applies to the set of densities $\mathcal{P}$. A first somewhat trivial but necessary observation is that $\mathcal{P}$ is nonempty; in the convolution kernel case this is already clear from the discussion before Proposition 4. Moreover, in the case of the Haar wavelet, it is not difficult to prove directly that if the density $f$ is such that, for some $x_0 \in \mathbb{R}$,

$$\frac{f(x_0 - v) - f(x_0)}{\mathrm{sign}(v)|v|^t}$$

has a nonnegative limit as $v \to 0$, $D$ say, then, for $k_0$ such that $x_0 \in (k_0/2^l, (k_0+1)/2^l]$ we have $|\beta_{lk_0}(f)| \geq 2^{-l(t+1/2)}(D-\varepsilon)(1+t)^{-1}(1-2^{-t})$ for every $\varepsilon > 0$ and $l = l(\varepsilon)$ large enough, which verifies the lower bound in Condition 3 by (4.45) below. In particular *every* differentiable density $f$ satisfies this condition for the Haar wavelet. For completeness we show that $\mathcal{P}$ is nonempty also for more general wavelet bases (see Section 4.4.2 below). In fact we shall see there that $\mathcal{P}$ is quite rich; small local modifications of arbitrary densities $f \in \mathcal{C}^t(\mathbb{R})$ are contained in $\mathcal{P}$.

To return to the interpretation of Proposition 4, let $\Sigma(t, L)$ be a ball in $\mathcal{C}^t(\mathbb{R})$, and define the subset of densities that are bounded away from zero on $F$,

$$\tilde{\Sigma}(t) := \Sigma(t, L) \cap \mathcal{D}(\eta, L, \delta, F).$$

It is natural to consider the *trace* (or relative) topology on $\tilde{\Sigma}(t)$ as a subset of the Banach space $\mathcal{C}^t(\mathbb{R})$. Proposition 4 then implies that the set of densities that are in $\tilde{\Sigma}(t)$ but *not* in $\mathcal{P}(\eta, r, b, b_1, b_2, \delta, F)$ for *any* $b_1, b_2$—so those *functions over which an adaptive sup-norm risk bound can be established but for which our adaptive confidence band is not necessarily valid*—is *nowhere dense* in the trace topology. If Theorem 1 is interpreted as a pointwise (in $f$) result, then these findings imply that there exists no (relatively) open set in $\tilde{\Sigma}(t)$ for which Condition 3 does not hold, and our adaptation result holds for "quasi-every" density in $\tilde{\Sigma}(t)$. Clearly, to obtain uniformity of the limit in Theorem 1, we have to fix a value of $b_1$, but *given* a (relatively) open set $\mathcal{O}$ in $\tilde{\Sigma}(t)$ we can always choose $b_1$ so small that $\mathcal{O}$ intersects $\mathcal{P}$ and hence is not contained in $\tilde{\Sigma}(t) \setminus \mathcal{P}$.

We finally remark that the minimax sup-norm risk over $\tilde{\Sigma}(t)$ is, for every $n$, the same as the one over $\tilde{\Sigma}(t) \cap \mathcal{P}$: This follows from the fact that $\mathcal{P} \cap \tilde{\Sigma}(t)$ is $\|\cdot\|_\infty$-dense (as $b_1 \searrow 0$) in $\tilde{\Sigma}(t)$, that the mapping $f \mapsto E_f \|T_n - f\|_\infty$ is continuous from $(\tilde{\Sigma}(t), \|\cdot\|_\infty)$ to $\mathbb{R}$ if $T_n$ is any measurable function of the sample which satisfies, for some fixed constant $c$, $\|T_n\|_\infty \leq c$ with probability one and that estimators $T_n$ that do not satisfy the last property for any $c$ can be neglected in the minimax risk.

**4. Proofs.** The proofs are organized into several parts. We start with some probabilistic results (Sections 4.1–4.3) that are central to verifying Condition 2. The statistically more relevant proofs are in Section 4.4 and depend on these probabilistic results only through Condition 2, so they can be read independently.

4.1. *A Gaussian reduction for maximal deviations of linear density estimators.* In what follows, given a metric space $(T, d)$ and $\varepsilon > 0$, the covering number $N(T, d, \varepsilon)$ denotes the smallest possible cardinality of any covering of $(T, d)$ by closed $d$-balls of radius at most $\varepsilon$. Its logarithm is referred to as the metric entropy of $(T, d)$. We also recall that a process $Y(t)$ on a metric space $(T, d)$ is said to be sample continuous if there exists a version of the process whose sample paths are all bounded and uniformly continuous. For Gaussian processes $Y(t)$, $t \in T$, unless specified otherwise, the distance $d$ is automatically taken to be the one provided by the process itself,

$$d(s,t) = (E(Y(t) - Y(s))^2)^{1/2}.$$

We let $K : \mathbb{R}^2 \mapsto \mathbb{R}$ be a measurable function satisfying:

(K1) $K$ is symmetric in its arguments, bounded, and for all $s \in \mathbb{R}$, $K(s,t)$ is right or left continuous in $t$ for every $t \in \mathbb{R}$;

(K2) $\sup_t \|K(t, \cdot)\|_V := \|K\|_V < \infty$ where $\|\cdot\|_V$ denotes the total variation norm on $\mathbb{R}$, $K(t, -\infty) = 0$ for all $t$;

(K3) there is a bounded, nonincreasing, exponentially decaying function $\Phi : \mathbb{R}^+ \cup \{0\} \mapsto \mathbb{R}^+ \cup \{0\}$ such that

$$|K(x,y)| \leq \Phi(|x-y|);$$

(K4) for all $\lambda \geq 1$, the covering numbers $N(\lambda[F_1, F_2], d, \varepsilon)$ of the intervals $[\lambda F_1, \lambda F_2]$ for the pseudo-distance $d(s,t) = (\int_\mathbb{R} (K(t,u) - K(s,u))^2 \, du)^{1/2}$ admit the bounds

$$N(\lambda[F_1, F_2], d, \varepsilon) \leq \frac{A' \lambda^{v_2}}{\varepsilon^{v_1}}$$

for some $A', v_i < \infty$ independent of $\varepsilon, \lambda$, and these bounds are valid for all positive $\varepsilon$ not exceeding the $d$ diameter of $[\lambda F_1, \lambda F_2]$, and (K5) there exist $\bar{A}, \bar{v}$ finite such that if $\mathcal{K} = \{K(2^j t, 2^j(\cdot)) : t \in \mathbb{R},\ j \in \mathbb{N} \cup \{0\}\}$ and if $\mathcal{Q}$ is the set of Borel probability measures on $\mathbb{R}$, then

$$\sup_{Q \in \mathcal{Q}} N(\mathcal{K}, L_2(Q), \varepsilon) \leq \left(\frac{\bar{A}}{\varepsilon}\right)^{\bar{v}} \tag{4.1}$$

for $0 < \varepsilon \leq \|K\|_\infty$.

Let $I = [a, b]$. Given a real sequence $j_n \to \infty$, define on $I$ the Gaussian processes

$$Y_n(t) = 2^{j_n/2} \int_\mathbb{R} K(2^{j_n} t, 2^{j_n} s)\, dW(s) = \int_\mathbb{R} K(2^{j_n} t, u)\, dW(u), \tag{4.2}$$

where $W$ is standard Brownian motion. It will often be convenient to rewrite $Y_n(t)$ as $Y_n(t) = Y(2^{j_n} t)$ where

$$Y(t) = \int_{-\infty}^{\infty} K(t, s)\, dW(s). \tag{4.3}$$

Note also that condition (K4) ensures that the processes $Y_n$ are sample continuous; for $u, v \in I$,

$$\begin{aligned} d_n^2(u, v) &:= E(Y_n(u) - Y_n(v))^2 \\ &= \int_\mathbb{R} (K(2^{j_n} u, s) - K(2^{j_n} v, s))^2\, ds \leq d^2(2^{j_n} u, 2^{j_n} v) \end{aligned} \tag{4.4}$$

so that $N(I, d_n, \varepsilon) \leq N(2^{j_n} I, d, \varepsilon)$ and it follows from condition (K4) that the square root of the metric entropy of $I$ with respect to the distance $d_n$ is integrable at zero, and hence the claim is an immediate consequence of Dudley's theorem [Theorem 2.6.1 in Dudley (1999)]. In particular, if we still denote a sample continuous version of $Y_n$ by $Y_n$, the norms $\|Y_n\|_I = \sup_{t \in I} |Y_n(t)|$ are proper random variables.

Let now

$$\mathcal{F}_n = \bigcup_{f \in \mathcal{D}} \mathcal{F}_n^f, \qquad \mathcal{F}_n^f = \{K(2^{j_n} t, 2^{j_n} \cdot)/\sqrt{f(t)} : t \in I\}. \tag{4.5}$$

Given $f \in \mathcal{D}$, let $X_i$ be i.i.d. with law $dP_f(t) := f(t)\, dt$, and let, as usual,

$$\nu_n^f = \frac{1}{\sqrt{n}} \sum_{i=1}^n (\delta_{X_i} - P_f), \qquad n \in \mathbb{N}$$

be the empirical processes based on the sequence $X_i$. Note that by the properties of $K$ and $f$, the supremum in $\|\nu_n^f\|_{\mathcal{F}_n^f}$ is countable and hence measurable.

The goal of this subsection is to prove the following proposition, in the spirit of Bickel and Rosenblatt (1973) and Giné, Koltchinskii and Sakhanenko (2004), and the proof adapts techniques from the last reference to the present situation. In what follows, $\Pr_f$ will still denote the product probability $P_f^{\mathbb{N}}$, but the symbol Pr will denote the probability measure determining the laws of all relevant other random variables (such as $Y_n$, and random variables constructed in the Gaussian coupling in the proof of Proposition 5 below).

PROPOSITION 5. *Let $I = [a,b]$, let $K$ be a function satisfying conditions (K1)–(K5) above, and let $j_n \to \infty$ as $n \to \infty$. Let $\{A_n\}$ and $\{B_n\}$ be numerical sequences such that $A_n \to \infty$ and*

(4.6) $$A_n = o\left(\frac{\sqrt{n}}{2^{j_n/2}\log n} \wedge 2^{j_n/2} \wedge \frac{2^{\alpha j_n}}{\sqrt{j_n}}\right)$$

*for some $0 < \alpha < 1$. Assume that there exists a random variable $Z$ with continuous distribution such that*

(4.7) $$\lim_{n\to\infty}\Pr\{A_n(\|Y_n\|_I - B_n) \leq x\} = \Pr\{Z \leq x\}, \qquad x \in \mathbb{R},$$

*where the processes $Y_n$ are defined by (4.2). Let $\mathcal{D}(\alpha, D, \delta, F)$ be as in (3.3) for the given $\alpha$, some $0 < D < \infty$, and admissible $\delta > 0, F = [F_1, F_2] \supset I$. Define, for each $f \in \mathcal{D}$, $\mathcal{F}_n^f$ as in (4.5), and let further $\nu_n^f$, $n \in \mathbb{N}$, be the empirical processes based on the variables $X_i$. Then, for all $x \in \mathbb{R}$,*

(4.8) $$\lim_{n\to\infty}\sup_{f\in\mathcal{D}}|\Pr_f\{A_n(2^{j_n/2}\|\nu_n^f\|_{\mathcal{F}_n^f} - B_n) \leq x\} - \Pr\{Z \leq x\}| = 0.$$

PROOF. *Step* 1: By Theorem 3 in Komlós, Major and Tusnády (1975), there is a probability space with a sequence $\{\xi_i\}$ of i.i.d. uniform on $[0,1]$ random variables and a sequence of Brownian motions $W_n$ defined on it such that, setting

$$\alpha_n(t) = \frac{1}{\sqrt{n}}\sum_{i=1}^n(\delta_{\xi_i}[0,t] - t)$$

and $W_n^\circ(t) = W_n(t) - tW_n(1)$, then,

(4.9) $$\Pr\left\{\|\alpha_n - W_n^\circ\|_{[0,1]} > \frac{x + C\log n}{\sqrt{n}}\right\} \leq \Lambda e^{-\theta x}, \qquad 0 \leq x < \infty, n \in \mathbb{N},$$

for some universal finite, positive constants $C$, $\Lambda$, $\theta$.

Define new random variables $\tilde{X}_i = F_f^{-1}(\xi_i)$, where $F_f^{-1}$ is the left continuous generalized inverse of the distribution function $F_f$ of $f$, right continuous at zero. For every $f \in \mathcal{D}$ the variables $\tilde{X}_i$ are i.i.d. with law $P_f$, and we denote by $\tilde{\nu}_n^f$ the associated empirical process. By (K2) and $f \geq \delta$ on $I$, the functions in $\mathcal{F}_n$ have total variation norm not exceeding $\|K\|_V/\sqrt{\delta}$, and since $F_f^{-1}$ is monotone, it follows that the same bound on the total variation norm (for functions on $[0,1]$) holds for all the functions in the classes

$$\tilde{\mathcal{F}}_n^f = \{h \circ F_f^{-1} : h \in \mathcal{F}_n^f\}, \qquad f \in \mathcal{D}, n \in \mathbb{N}.$$

Moreover, if $g$ is nonincreasing on $[0,1]$ with $g(0) = 1$ and $g(1) = 0$, then $g$ is the pointwise nondecreasing limit—and by dominated convergence, also the limit in $L_2([0,1])$—of convex combinations of indicators $I_{[0,t]}$, $0 \leq t \leq 1$. So by (K2) both $\alpha_n$ and $W_n$ extend from indicator functions $I_{[0,t]}$ to functions in $\tilde{\mathcal{F}}_n^f$ by linearity and continuity [see, e.g., Dudley (1985), Theorems 5.1–5.3], and so does $W_n^\circ$. We conclude that, for all $f \in \mathcal{D}$,

$$\|\alpha_n - W_n^\circ\|_{\tilde{\mathcal{F}}_n^f} \leq \|K\|_V \delta^{-1/2} \|\alpha_n - W_n^\circ\|_{[0,1]},$$

and, writing $G_{n,f}^\circ(g) = W_n^\circ(g \circ F_f^{-1})$ for $g \in \mathcal{F}_n^f$, that $E(G_{n,f}^\circ(g)G_{n,f}^\circ(\bar{g})) = P_f(g\bar{g}) - (P_f g)(P_f \bar{g})$, i.e., $G_{n,f}^\circ$ is a (sample continuous) version of the $P_f$-Brownian bridge. Since, furthermore,

$$\alpha_n(g \circ F_f^{-1}) = \tilde{\nu}_n^f(g)$$

by construction, (4.9) gives

$$\sup_{f \in \mathcal{D}} \Pr\left\{ \|\tilde{\nu}_n^f - G_{n,f}^\circ\|_{\mathcal{F}_n^f} > \frac{\|K\|_V \delta^{-1/2}(x + C\log n)}{\sqrt{n}} \right\}$$

$$= \sup_{f \in \mathcal{D}} \Pr\left\{ \|\alpha_n - W_n^\circ\|_{\tilde{\mathcal{F}}_n^f} > \frac{\|K\|_V \delta^{-1/2}(x + C\log n)}{\sqrt{n}} \right\}$$

$$\leq \Pr\left\{ \|\alpha_n - W_n^\circ\|_{[0,1]} > \frac{x + C\log n}{\sqrt{n}} \right\} \leq \Lambda e^{-\theta x}$$

for all $x \geq 0$ and $n \in \mathbb{N}$. Taking $x = (C' - C)\log n$ for some $C' > C$ in this inequality, we have

$$\sup_{f \in \mathcal{D}} \Pr\left\{ A_n 2^{j_n/2} \|\tilde{\nu}_n^f - G_{n,f}^\circ\|_{\mathcal{F}_n^f} > \frac{\|K\|_V \delta^{-1/2} C' A_n \log n}{\sqrt{n 2^{-j_n}}} \right\}$$

(4.10)
$$\leq \frac{\Lambda}{n^{(C'-C)\theta}}.$$

In particular, if

(4.11)
$$A_n = o\left(\frac{\sqrt{n}}{2^{j_n/2} \log n}\right),$$

then (4.10) implies that there exists a sequence $\varepsilon_n'' \to 0$ such that

(4.12)
$$\limsup_n \sup_{f \in \mathcal{D}} \Pr\{A_n 2^{j_n/2} \|\tilde{\nu}_n^f - G_{n,f}^\circ\|_{\mathcal{F}_n^f} \geq \varepsilon_n''\} = 0.$$

Consider next the processes $G_{n,f}(g) = W_n(g \circ F_f^{-1})$, $g \in \mathcal{F}_n^f$ which are sample continuous versions of the $P_f$-Brownian motion $W_{P_f}$ since

$$E(W_n(g \circ F_f^{-1})W_n(\bar{g} \circ F_f^{-1})) = \int_0^1 g \circ F_f^{-1}(x)\bar{g} \circ F_f^{-1}(x)\, dx$$

$$= \int_\mathbb{R} g(y)\bar{g}(y)f(y)\, dy.$$

Since $W_n^\circ(g \circ F^{-1}) = W_n(g \circ F^{-1}) - (\int_0^1 g \circ F^{-1}(t)\, dt)W_n(1)$, and, since, by (K3),

$$(4.13) \qquad \sup_{g \in \mathcal{F}_n^f} \left| \int_0^1 g \circ F^{-1}(t)\, dt \right| = \sup_{g \in \mathcal{F}_n^f} |P_f g| \leq \delta^{-1/2}\|\Phi\|_1 2^{-j_n},$$

it follows that if

$$(4.14) \qquad A_n = o(2^{j_n/2}),$$

then we can replace $G_n^\circ$ by $G_n$ in (4.12); that is, there exists $\varepsilon_n' \to 0$ such that

$$(4.15) \qquad \limsup_n \sup_{f \in \mathcal{D}} \Pr\{A_n 2^{j_n/2}\|\tilde{\nu}_n^f - G_{n,f}\|_{\mathcal{F}_n^f} \geq \varepsilon_n'\} = 0.$$

[Note that by the results of Dudley (1985) alluded to above, for all $n$ and $f$, the process $W_{P_f}(g)$, $g \in \mathcal{F}_n^f$, is sample continuous (hence sample bounded).]

*Step* 2: To compare $G_{n,f}$ on $\mathcal{F}_n^f$ with $Y_n$ we must couple in the right way sample continuous versions of both processes. Since the functions in $\mathcal{F}_n^f$ are parametrized by $t \in I$, we will write (in slight abuse of notation) $G_{n,f}(t)$, $t \in I$, for $G_{n,f}(g_t)$, $g_t(\cdot) = K(2^{j_n}t, 2^{j_n}\cdot)/\sqrt{f(t)} \in \mathcal{F}_n^f$. First we observe that the process,

$$W(K(2^{j_n}t, 2^{j_n}\cdot)\sqrt{f(\cdot)/f(t)}), \qquad t \in I,$$

where $W$ is Brownian motion acting on functions as described in step 1, is a version of $G_{n,f}$ (both processes have the same covariance). Next we observe that the set $\mathcal{G}_n$ defined by

$$\mathcal{G}_n = \{2^{j_n/2}K(2^{j_n}t, 2^{j_n}\cdot), K(2^{j_n}t, 2^{j_n}\cdot)\sqrt{f(\cdot)/f(t)} : t \in I\}$$

is a GC subset of $L_2(\mathbb{R}, \lambda)$ where $\lambda$ is Lebesgue measure; this follows from the entropy bounds (K4) and (K5) and the results in Dudley (1999), Sections 2.5 and 2.6, in particular Theorem 2.6.1. Hence, the restriction to $\mathcal{G}_n$ of the isonormal process of $L_2(\mathbb{R}, \lambda)$, which we write as $g \mapsto \int_\mathbb{R} g\, dW = W(g)$, admits a version with bounded uniformly continuous sample paths [for the $L_2(\mathbb{R}, \lambda)$ distance]. We call $\tilde{G}_{n,f}(t)$ and $\tilde{Y}_n(t)$ the restrictions of this process to the sets $\{K(2^{j_n}t, 2^{j_n}\cdot)\sqrt{f(\cdot)/f(t)}: t \in I\}$ and $\{2^{j_n/2}K(2^{j_n}t, 2^{j_n}\cdot): t \in I\}$, respectively. They are versions of $G_{n,f}$ and $Y_n$, respectively, and, as we see next, we can control the supnorm of their difference. Set

$$Z_{n,f}(t) = 2^{j_n/2}\tilde{G}_{n,f}(t) - \tilde{Y}_n(t)$$

$$= 2^{j_n/2}\int_{\mathbb{R}} K(2^{j_n}t, 2^{j_n}s)\left(\sqrt{\frac{f(s)}{f(t)}} - 1\right) dW(s), \quad t \in I.$$

We have for $u, v \in I$,

$$d_{Z_{n,f}}(u,v) := (E(Z_{n,f}(u) - Z_{n,f}(v))^2)^{1/2}$$

$$\leq \delta^{-1/2}\|K(2^{j_n}u, \cdot) - K(2^{j_n}v, \cdot)\|_{L_2(P_f)}$$

$$+ \|K\|_\infty |f^{-1/2}(u) - f^{-1/2}(v)| + d_n(u,v),$$

where $d_n(u,v) = d(2^{j_n}u, 2^{j_n}v)$ [cf. (4.4)], and, using (K4), (K5) and that $\delta \leq f(u) \leq D$ for all $u \in I$, it is easy to show that the covering numbers of $I$ for this distance are bounded by

(4.16) $$N(I, d_{Z_{n,f}}, \varepsilon) \leq B2^{v_3 j_n}/\varepsilon^{v_4}$$

for every (small) $\varepsilon > 0$ and constants $B, v_3, v_4$ independent of $n$ and $f$. Note next that by Remark 1, for every $f \in \mathcal{D}$ there exists $L = L(\mathcal{D})$ and $c = c(\mathcal{D})$ such that $\|f\|_\infty \leq L$ and the $\alpha$-Hölder constant of $f$ is at most $c$. Hence, for $t \in I$,

$$(\sqrt{f(t - 2^{-j_n}u)} - \sqrt{f(t)})^2 \leq LI(|u| > \delta 2^{j_n}) + c^2 2^{-2\alpha j_n} I(|u| \leq \delta 2^{j_n}),$$

and we obtain for all $t \in I$,

(4.17)
$$E(Z_{n,f}(t))^2 = 2^{j_n} \int_{\mathbb{R}} K^2(2^{j_n}t, 2^{j_n}s)\left(\sqrt{\frac{f(s)}{f(t)}} - 1\right)^2 ds$$

$$\leq 2^{j_n}\delta^{-1}\int_{\mathbb{R}} K^2(2^{j_n}t, 2^{j_n}s)(\sqrt{f(s)} - \sqrt{f(t)})^2 ds$$

$$\leq \delta^{-1}\int_{\mathbb{R}} \Phi^2(u)(\sqrt{f(t - 2^{-j_n}u)} - \sqrt{f(t)})^2 du$$

$$\leq \delta^{-1}L\|\Phi\|_1 \Phi(\delta 2^{j_n}) + \delta^{-1}\|\Phi\|_2^2 c^2 2^{-2\alpha j_n} \leq D_1^2 2^{-2\alpha j_n},$$

where $D_1$ is a constant that does not depend on $n$ nor $f$. That is, the diameter of $I$ for the $L_2$-distance induced by the process $Z_{n,f}$ is at most $2D_1 2^{-\alpha j_n}$. Hence, Dudley's entropy bound in expectation form [e.g., de la

Peña and Giné (1999), Corollary 5.1.6], (4.16) and (4.17) give

$$E \sup_{t \in I} |2^{j_n/2} \tilde{G}_{n,f}(t) - \tilde{Y}_n(t)| \lesssim D_1 2^{-\alpha j_n} + \int_0^{D_1 2^{-\alpha j_n}} \sqrt{\log \frac{B 2^{v_3 j_n}}{\varepsilon^{v_4}}} \, d\varepsilon$$

$$\lesssim \sqrt{j_n} 2^{-\alpha j_n}$$

with unspecified constants independent of $f \in \mathcal{D}$ and $n$. So, if, besides (4.11) and (4.14), the sequence $\{A_n\}$ satisfies

$$A_n = o(2^{\alpha j_n}/\sqrt{j_n})$$

[hence, if $\{A_n\}$ satisfies (4.6)], then there exists $\varepsilon_n \to 0$ such that

(4.18) $$\lim_{n \to \infty} \sup_{f \in \mathcal{D}} \Pr\{A_n \|2^{j_n/2} \tilde{G}_{n,f} - \tilde{Y}_n\|_I \geq \varepsilon_n\} = 0.$$

*Step* 3: We finally combine the bounds obtained. Clearly $\|\tilde{G}_{n,f}\|_I$ has the same probability law as $\|G_{n,f}\|_{\mathcal{F}_n^f}$, and likewise $\|\tilde{Y}_n\|_I$ has the same law as $\|Y_n\|_I$. Therefore, under the hypotheses of the proposition, we have, for all $f \in \mathcal{D}$ and $x_n \to x$, $x \in \mathbb{R}$,

$$[\Pr\{A_n(\|Y_n\|_I - B_n) \leq x_n - \varepsilon_n\} - \Pr\{Z \leq x\}]$$
$$- \sup_{f \in \mathcal{D}} \Pr\{A_n \|2^{j_n/2} \tilde{G}_{n,f} - \tilde{Y}_n\|_I > \varepsilon_n\}$$
$$\leq \Pr\{A_n(2^{j_n/2}\|G_{n,f}\|_{\mathcal{F}_n^f} - B_n) \leq x_n\} - \Pr\{Z \leq x\}$$
$$\leq [\Pr\{A_n(\|Y_n\|_I - B_n) \leq x_n + \varepsilon_n\} - \Pr\{Z \leq x\}]$$
$$+ \sup_{f \in \mathcal{D}} \Pr\{A_n \|2^{j_n/2} \tilde{G}_{n,f} - \tilde{Y}_n\|_I > \varepsilon_n\}.$$

The leftmost and rightmost sides of this inequality do not depend on $f \in \mathcal{D}$ and tend to zero by (4.7), the continuity of the probability law of $Z$ and (4.18). Thus we have

(4.19) $$\lim_{n \to \infty} \sup_{f \in \mathcal{F}_n^f} |\Pr\{A_n(2^{j_n/2}\|G_{n,f}\|_{\mathcal{F}_n^f} - B_n) \leq x_n\} - \Pr\{Z \leq x\}| = 0$$

for any sequence $x_n \to x$, any $x \in \mathbb{R}$. Similarly, since the random variables $\|\tilde{\nu}_n^f\|_{\mathcal{F}_n^f}$ and $\|\nu_n^f\|_{\mathcal{F}_n^f}$ have the same law, we have, for any $x \in \mathbb{R}$,

$$[\Pr\{A_n(2^{j_n/2}\|G_{n,f}\|_{\mathcal{F}_n^f} - B_n) \leq x - \varepsilon_n'\} - \Pr\{Z \leq x\}]$$
$$- \sup_{f \in \mathcal{D}} \Pr\{A_n 2^{j_n/2} \|\tilde{\nu}_n^f - \tilde{G}_{n,f}\|_{\mathcal{F}_n^f} > \varepsilon_n'\}$$
$$\leq \Pr_f\{A_n(2^{j_n/2}\|\nu_n^f\|_{\mathcal{F}_n^f} - B_n) \leq x\} - \Pr\{Z \leq x\}$$

$$\leq [\Pr\{A_n(2^{j_n/2}\|G_{n,f}\|_{\mathcal{F}_n^f} - B_n) \leq x + \varepsilon_n'\} - \Pr\{Z \leq x\}]$$
$$+ \sup_{f \in \mathcal{D}} \Pr\{A_n 2^{j_n/2}\|\tilde{\nu}_n^f - \tilde{G}_{n,f}\|_{\mathcal{F}_n^f} > \varepsilon_n'\}$$

which, by (4.15) and by (4.19) with $x_n = x \pm \varepsilon_n'$, gives (4.8). □

Condition (K3) is only used in the equation above (4.13) and in (4.17); therefore it can be relaxed to the following: there is $\Phi$ measurable, bounded and satisfying that, for some $y_0$ and $\eta > 0$ and all $y > y_0$, $\sup_{x \geq y} \Phi(x) \leq y^{-1-\eta}$, such that $|K(x,y)|$ is dominated by $\Phi(|x-y|)$.

4.2. *Examples and some first limit theorems for suprema of certain Gaussian processes.*

4.2.1. *Haar wavelets.* The projection kernel corresponding to the Haar wavelet is

(4.20)  $$K(x,y) = \sum_{k \in \mathbb{Z}} I_{[0,1)}(x-k)I_{[0,1)}(y-k) = I([x] = [y]).$$

It obviously satisfies conditions (K1)–(K3) ($\|K(t,\cdot)\|_V = 2$, $\Phi(|u|) = I(|u| \leq 1)$). Moreover, $d^2(x,y) = \int_\mathbb{R}(K(x,u) - K(y,u))^2 du = 0$ if $[x] = [y]$ and 2 otherwise so that

$$N(\lambda[F_1, F_2], d, \varepsilon) \leq N(\lambda[F_1, F_2], d, 0) \leq \lambda(F_2 - F_1) + 2 \leq \frac{2\lambda(F_2 - F_1) + 4}{\varepsilon}$$

for $0 < \varepsilon < 2$ (note that 2 is an upper bound for the $d$-diameter of any set of real numbers) so that (K4) holds. Condition (K5) follows from Lemma 2 in Giné and Nickl (2009b).

So Proposition 5 applies and we are led to consider the process [see (4.3)],

$$Y(t) = \sum_{k \in \mathbb{Z}} I(t \in [k, k+1)) \int_k^{k+1} dW(s) = \sum_{k \in \mathbb{Z}} I(t \in [k, k+1))g_k,$$

where $g_k$ are i.i.d. $N(0,1)$, and therefore, taking $I = [0,1]$,

$$\sup_{0 \leq t \leq 1} |Y_n(t)| = \sup_{0 \leq u \leq 2^{j_n}} |Y(u)| = \max_{0 \leq k \leq 2^{j_n}} |g_k|.$$

Now Theorem 1.5.3 (and Theorem 1.8.3) in Leadbetter et al. (1983) gives

$$\Pr\left\{A_n\left(\sup_{0 \leq t \leq 1} |Y_n(t)| - B_n\right) \leq x\right\} \to e^{-e^{-x}} \qquad \text{for all } x \in \mathbb{R},$$

where $A_n = A(j_n)$, $B_n = B(j_n)$ and

(4.21)  $$A(l) = [2(\log 2)l]^{1/2}, \qquad B(l) = A(l) - \frac{\log l + \log(\pi \log 2)}{2A(l)}.$$

Combining this with Proposition 5 we have, recalling the set $\mathcal{D}$ from (3.3).

PROPOSITION 6. *Let $\mathcal{D} = \mathcal{D}(\alpha, D, \delta, F)$ for some $0 < \alpha \leq 1$, $0 < D < \infty$, and where $\delta, F$ are admissible. If $j_n \to \infty$ as $n \to \infty$ satisfies $j_n 2^{j_n} = o(n/(\log n)^2)$ and if $f_n := f_n(\cdot, j_n)$ is the Haar wavelet estimator from (2.1) with $\phi = 1_{[0,1)}$, then*

$$\sup_{f \in \mathcal{D}} \left| \Pr_f \left\{ A_n \left( \sqrt{n 2^{-j_n}} \left\| \frac{f_n - E f_n}{\sqrt{f}} \right\|_{[0,1]} - B_n \right) \leq x \right\} - e^{-e^{-x}} \right| \to 0$$

*for all $x \in \mathbb{R}$*

*as $n \to \infty$ where $A_n$ and $B_n$ are as above (4.21).*

4.2.2. *Convolution kernels.* If $K$ is a real-valued function with bounded support, and is symmetric and Lipschitz continuous, then the kernel $K(x, y) := K(x - y)$ satisfies conditions (K1)–(K4) with $\Phi = K$ and $d(s, t)$ proportional to $|s - t|$. Condition (K5) is proved in Nolan and Pollard (1987). [These are not the only convolution kernels satisfying (K1)–(K5) and, for instance, the Gaussian kernel also satisfies theses conditions.]

Assume in what follows that $K$ is bounded, symmetric, supported by $[-1, 1]$ and twice continuously differentiable on $\mathbb{R}$. Writing $Y_n(t) = Y(2^{j_n} t)$ with $Y$ as in (4.3) we have

$$\sup_{t \in [0,1]} |Y_n(t)| = \sup_{0 \leq t \leq 2^{j_n}} |Y(t)|.$$

In this case, $Y(t) = \int_{\mathbb{R}} K(t - s)\, dW(s)$ is a stationary Gaussian process with covariance

$$r(t) := E(Y(t) Y(0)) = \int_{\mathbb{R}} K(t + u) K(u)\, du = \|K\|_2^2 - C t^2 + o(t^2),$$

where $C = -2^{-1} \int_{\mathbb{R}} K(u) K''(u)\, du > 0$ (by integration by parts), and $r(t) = 0$ for $|t| > 2$. Set $\tilde{Y} = Y/\|K\|_2$ and $\tilde{C} = C/\|K\|_2^2$. We apply Theorem 8.2.7 in Leadbetter et al. (1983), in its version for absolute values (Corollary 11.1.6 in the same reference) with

$$(4.22) \qquad B(l) = \sqrt{2(\log 2) l} + \frac{\log \sqrt{2 \tilde{C}} - \log \pi}{\sqrt{2(\log 2) l}}.$$

One has

$$\lim_{n \to \infty} \Pr\left\{ \sqrt{2(\log 2) j_n} \left( \sup_{0 \leq t \leq 2^{j_n}} |Y(t)|/\|K\|_2 - B(j_n) \right) \leq x \right\} \to e^{-e^{-x}}, \qquad x \in \mathbb{R}$$

which, combined with Proposition 5, yields the following:

PROPOSITION 7. *If $K:\mathbb{R} \mapsto \mathbb{R}$ is bounded, symmetric, supported by $[-1,1]$ and twice continuously differentiable, $\mathcal{D}$ and $j_n$ are as in Corollary 6, $B(l)$ is as in (4.22) and if $f_n := f_n(y, j_n)$ is the kernel estimator from (2.1), then, as $n \to \infty$,*

$$\sup_{f \in \mathcal{D}} \left| \Pr_f \left\{ \sqrt{2(\log 2) j_n} \left( \sqrt{n 2^{-j_n}} \left\| \frac{f_n - E f_n}{\|K\|_2 \sqrt{f}} \right\|_{[0,1]} - B(j_n) \right) \leq x \right\} - e^{-e^{-x}} \right| \to 0$$

*for all $x \in \mathbb{R}$.*

4.2.3. *General wavelet bases.* Let $K(x, y) = \sum_k \phi(x - k) \phi(y - k)$ be a general wavelet projection kernel with scaling function $\phi$. Assuming the conditions of Proposition 5 are verified for the moment (see below for examples), we are led to consider the distributions of maxima over increasing intervals of the process

$$Y(t) = \int_\mathbb{R} K(t, u) \, dW(u) = \sum_{k \in \mathbb{Z}} \phi(t - k) \int_{-\infty}^\infty \phi(s - k) \, dW(s),$$

where $W$ is Brownian motion. Since the functions $\phi(\cdot - k)$ are orthonormal in $L_2(\mathbb{R})$, we can write this process as

(4.23) $$Y(t) = \sum_{k \in \mathbb{Z}} \phi(t - k) g_{-k} = \sum_{k \in \mathbb{Z}} \phi(t + k) g_k,$$

where $g_k = \int_\mathbb{R} \phi(u + k) \, dW(u)$ is a sequence of i.i.d. standard normal random variables.

The process $Y(t)$ is *not* stationary in general. However, the covariance of this process, which is given by

(4.24) $$r(t, t + v) := EY(t) Y(t + v) = \sum_k \phi(t - k) \phi(t + v - k)$$

is, for each $v \in \mathbb{R}$, periodic with period one, and so is its variance function,

(4.25) $$\sigma^2(t) := EY^2(t) = E \left( \sum_k \phi(t - k) g_k \right)^2 = \sum_k \phi^2(t - k).$$

Processes $Y(t), t \in \mathbb{R}$, whose covariance function $t \mapsto r(t, t + v)$ is periodic in $t$ for every $v \in \mathbb{R}$, with a period independent of $v$, are called *cyclostationary*. So, although $Y$ is not stationary, it is cyclostationary with period one.

While the asymptotic distributional theory for maxima of nonstationary Gaussian processes is not as complete as for stationary processes, there are in particular some results for "cyclostationary processes." These results involve a careful analysis of the variance and covariance functions, and this requires a case-by-case treatment of the wavelet basis in question. For Battle–Lemarié wavelets we can carry this through and will prove the relevant limit theorem for their associated Gaussian processes $Y$ in the next section.

Particularly interesting are the compactly supported wavelets, such as the Daubechies family. The Daubechies scaling functions $\phi = \phi_N$ [see Daubechies (1992) and also Chapter 7 in Härdle et al. (1998)] have compact support and are differentiable for $N \geq 3$. The associated projection kernel $K$ obviously satisfies (K1) and (K2); it satisfies (K3), as shown, for example, in Härdle et al. (1998), Lemmas 8.5 and 8.6; (K5) is proved in Lemma 2 in Giné and Nickl (2009b); and regarding (K4) we note that, by orthonormality,

$$
\begin{aligned}
\int_{\mathbb{R}} (K(x,u) - K(y,u))^2 \, du &= \sum_k (\phi(x-k) - \phi(y-k))^2 \\
&\leq 2 \left\| \sum_k |\phi(\cdot - k)| \right\|_\infty \|\phi'\|_\infty |x-y|,
\end{aligned}
\tag{4.26}
$$

a distance which satisfies (K4), since $\|\sum_k |\phi(\cdot - k)|\|_\infty < \infty$ (e.g., Härdle et al. (1998), Lemma 8.5). So the Gaussian reduction in Proposition 5 applies to Daubechies-wavelets for $N \geq 3$, and therefore it remains to prove a limit theorem for the process $Y(t) = \int_{\mathbb{R}} K(t,u) \, dW(u)$. The covariance function of $Y$ in the case of Daubechies wavelets seems difficult to analyze and we do not know how to derive an exact distributional limit theorem for $\max_{0 \leq t \leq 2^{jn}} |Y(t)|$ in this case. We conjecture that these limits exist and are similar to the Battle–Lemarié cases considered in the next subsection.

### 4.3. The limit theorem for Battle–Lemarié wavelets.

#### 4.3.1. Battle–Lemarié wavelets. Let

$$N_r(x) = I_{[0,1)} * \cdots^{(r)} * I_{[0,1)}(x)$$

be the $B$-spline of order $r$, $r \in \mathbb{N}$. The scaling function $\phi_r$ that generates the Battle–Lemarié wavelet basis admits a unique representation,

$$\phi_r(x) = \sum_{k \in \mathbb{Z}} a_k^{(r)} N_r(x - k), \tag{4.27}$$

where the sequence of coefficients $a_k^{(r)}$ has exponential decay as $|k| \to \infty$ for all $r$ [see Daubechies (1992), Corollary 4.5.2, where one can also find an explicit definition of these coefficients]. For $r = 1$, $\phi_1$ is the Haar scaling function which has already been considered. For $r > 1$, the function $\phi_r$ is Lipschitz and, although it does not have bounded support, it decreases exponentially. In fact, as is easily seen, $\phi$ is a $(r-1)$-regular wavelet basis (cf. prior to Definition 1).

As a first step we verify the conditions of Proposition 5 for the kernel $K(x,y)$. Condition (K1) is obvious; condition (K3) holds with an exponentially decreasing $\Phi$ [using Lemma 8.6 in Härdle et al. (1998)]. Condition (K2) is, for example, contained in the proof of Lemma 2 in Giné and Nickl (2010) which itself verifies (K5). Regarding (K4) we can argue directly as in (4.26) for $r > 2$, in which case $\phi_r$ is differentiable and has uniformly bounded derivatives. In case $r = 2$ a similar argument works, using that $\phi_2$ is uniformly bounded and Lipschitz on $[k, k+1]$ for each $k \in \mathbb{Z}$.

Consequently Proposition 5 applies if we can derive a limit of type (4.7) for the process,

$$(4.28) \qquad Y^{(r)}(t) = \int_{\mathbb{R}} K(t,u)\,dW(u) = \sum_k \phi_r(x+k)g_k,$$

with variance and covariance as in Section 4.2.3 and where the $g_k$'s are i.i.d. standard normal.

Since for $r = 2$ the Battle–Lemarié scaling function $\phi_2$ is not differentiable, we will have to treat the cases $r = 2$ and $r > 2$ separately. In both cases we rely on the following elementary (but somewhat cumbersome) key lemma on the maxima of the variance function (4.25). It will be proved in Section 4.4.3.

LEMMA 1. *Let $\phi_r$ be the scaling function for the Battle–Lemarié wavelet of order $r$, $r = 2, 3, 4$ and let*

$$\sigma_r^2(t) = \sum_{k \in \mathbb{Z}} \phi_r^2(t-k), \qquad t \in \mathbb{R}.$$

*Then, $\sigma_3^2(t)$ attains its absolute maximum only at the points $t = l + 1/2$, $l \in \mathbb{Z}$, and for $r = 2, 4$, $\sigma_r^2(t)$ attains its absolute maximum only at the points $t = l \in \mathbb{Z}$. The values are, for every $l \in \mathbb{N}$,*

$$\sigma_2^2 := \sigma_2^2(l) = \sum_k \phi_2^2(k) = \sum_k (a_k^{(2)})^2,$$

$$\sigma_4^2 := \sigma_4^2(l) = \sum_k \phi_4^2(k) = \frac{1}{36}\sum_k (a_{k-1}^{(4)} + 4a_k^{(4)} + a_{k+1}^{(4)})^2$$

*and*

$$\sigma_3^2 := \sigma_3^2(l - 1/2) = \frac{M}{32} + \frac{1}{4}\sum_k (a_k^{(3)} + a_{k-1}^{(3)})^2,$$

*where*

$$M = \sum_k (a_k^{(3)} - a_{k-1}^{(3)})^2 - \sum_k (a_k^{(3)} - a_{k-1}^{(3)})(a_{k-1}^{(3)} - a_{k-2}^{(3)})$$

*and $a_k^{(r)}$, $k \in \mathbb{Z}$, $r = 2, 3, 4$, are the coefficients in (4.27).*

4.3.2. *The limit theorem for the piecewise linear Battle–Lemarié wavelet.*
We first treat the piecewise linear case $r = 2$ where the structure of the
process $Y^{(2)}$ allows for a direct reduction to the maximum of a stationary
Gaussian sequence. Since $Y^{(2)}(t)$ is piecewise linear, it can be written as
[cf. (4.47) below]

$$Y^{(2)}(t) = \sum_{k \in \mathbb{Z}} \phi_2(t+k) g_k = \sum_{k \in \mathbb{Z}} \phi_2(t+k-l) g_{k-l}$$

(4.29)
$$= (t-l) \sum_k a_k^{(2)} g_{k-l} + (1-(t-l)) \sum_k a_{k-1}^{(2)} g_{k-l},$$

$$l \leq t \leq l+1, l \in \mathbb{Z},$$

where the variables $g_k$ are i.i.d. standard normal. Let $X(t) := Y^{(2)}(t)/\sigma_2$
be the same process normalized by the maximum of the variance function
$\sigma_2(t)$ which is attained at $t = l \in \mathbb{Z}$ and given in Lemma 1. We shall write
$a_k = a_k^{(2)}$ and $\sigma = \sigma_2$ throughout this subsection, and we recall that the
coefficients $a_k$ decrease exponentially. On any interval $[l, l+1]$, $l \in \mathbb{Z}$, $X$ has
the form $X(t) = (t-l)G_1 + G_2$ and, obviously, the absolute maximum of
$|(t-l)G_1 + G_2|$ over $t \in [l, l+1]$ is attained either at $t = l$ or at $t = l+1$.
Therefore,

$$\sup_{t \in [0, 2^{jn}]} |X(t)| = \max_{0 \leq l \leq 2^{jn}, l \in \mathbb{Z}} |X(l)|.$$

Considering the process $X$ indexed only by integers, we see that, for all
$l, m \in \mathbb{Z}$,

$$\sigma^2 E(X(l)X(l+m)) = E\left(\sum_k a_{k-1} g_{k-l} \sum_{k'} a_{k'-1} g_{k'-l-m}\right)$$

$$= \sum_k a_{k+l-1} a_{k+l+m-1} = \sum_k a_k a_{k+m};$$

that is, the sequence $X(l)$, $l \in \mathbb{Z}$, is *stationary* with covariance

$$r(m) = \sum_k a_k a_{k+m} \Big/ \sum_k a_k^2.$$

Using the exponential decay of the $a_k$'s in the bound

$$\sum_k |a_k a_{k+m}| \leq 2 \Big(\sup_k |a_k|\Big) \sum_{|k| \geq m/2} |a_k|,$$

one sees $(\log m) r(m) \to 0$ as $m \to \infty$ (Berman's condition). Then we can apply the usual result about weak convergence of the maximum of a stationary Gaussian sequence satisfying Berman's condition [e.g., Theorem 4.3.3 in Leadbetter et al. (1983)], together with the asymptotic independence of maximum and minimum of such sequences [e.g., Davis (1979), page 459f]. The outcome is that the limit theorem for $\max_{0\leq l\leq 2^{j_n}, l\in \mathbb{Z}}|X(l)|$ is the same as if the $X(l)$ were independent; with $A_n, B_n$ as above (4.21),

$$\Pr\Big\{A_n\Big(\sup_{0\leq t\leq 2^{j_n}}|Y^{(2)}(t)|/\sigma - B_n\Big)\leq x\Big\} \to e^{-e^{-x}}, \qquad x\in\mathbb{R},$$

which, combined with Proposition 5 (which applies by observations in Section 4.3.1), gives the following proposition:

PROPOSITION 8. *Let $\mathcal{D}$, $j_n$ be as in Corollary 6. If $f_n^{(2)} := f_n(y, j_n)$ is the wavelet density estimator from (2.1) with Battle–Lemarié scaling function $\phi_2$, if $A_n = A(j_n)$ and $B_n = B(j_n)$ are as in (4.21) and if $\sigma^2 = \sum_k a_k^2$, then, as $n\to\infty$,*

$$\sup_{f\in\mathcal{D}}\left|\Pr_f\left\{A_n\left(\sqrt{n2^{-j_n}}\left\|\frac{f_n^{(2)} - Ef_n^{(2)}}{\sigma\sqrt{f}}\right\|_{[0,1]} - B_n\right)\leq x\right\} - e^{-e^{-x}}\right| \to 0$$

*for all $x\in\mathbb{R}$.*

4.3.3. *The limit theorem for smooth Battle–Lemarié wavelets.* In the case of $r>2$ we have to deal with maxima of nonstationary Gaussian processes. The following theorem is an adaptation to mean square differentiable cyclostationary processes [cf. after (4.25) above] of a theorem of Piterbarg and Seleznjev (1994) [see also Konstant and Piterbarg (1993), Hüsler (1999), Hüsler, Piterbarg and Seleznjev (2003)].

THEOREM 2. *Let $X(t)$, $t\in\mathbb{R}$, be a cyclostationary, centered Gaussian process with period 1, variance $\sigma_X(t)$ and covariance $r_X(s,t)$. Assume:*

(1) $X(t)$ *is mean square differentiable and a.s. continuous;*
(2) $r_X(s,t) = \sigma_X(s)\sigma_X(t)$ *only at $s=t$;*
(3) $\sup_{t\in[0,1]}\sigma_X^2(t) = \sigma_X^2(t_0) = 1$ *for a unique $t_0 \in (0,1)$, $\sigma_X^2(t)$ is twice continuously differentiable at $t_0$, $\sigma'_X(t_0) = 0$, $\sigma''_X(t_0) < 0$ and $E(X'(t_0))^2 > 0$;*
(4) $(\log v)\sup_{s,t:|s-t|\geq v}|r_X(s,t)|\to 0$ *as $v\to\infty$.*

*Then, for all $x\in\mathbb{R}$,*

$$\lim_{T\to\infty}\Pr\Big\{\sup_{t\in[0,T]}|X(t)|\leq \frac{x}{a_T} + b_T\Big\} = e^{-e^{-x}},$$

*where*

$$a_T = \sqrt{2\log T} \quad \text{and}$$
$$b_T = a_T - \frac{\log\log T + \log\pi - \log\sqrt{1 - E(X'(t_0))^2/\sigma''_X(t_0)}}{2a_T}.$$

PROOF. Under the hypotheses (1)–(3), the correlation $R_X(s,t) = r_X(s,t)/(\sigma_X(s)\sigma_X(t))$ admits the following development near $(t_0, t_0)$:

$$R_X(s,t) = 1 - \left(\frac{E(X'(t_0))^2}{2} + D(t,s)\right)(t-s)^2 + o((t-s)^2),$$

where $D(t,s)$ is continuous at $(t_0, t_0)$ and satisfies $D(t_0, t_0) = 0$ [see, e.g., the proof of Corollary 2.1 in Konstant and Piterbarg (1993)]. Because of this all the hypotheses of Theorem 1 in Hüsler (1999), or of Theorem 1 in Piterbarg and Seleznjev (1994) are satisfied with $2a = -\sigma''(t_0)$ and $2b = E(X'(t_0))^2$. Therefore, as $T \to \infty$,

$$\Pr\left\{\sup_{t \in [0,T]} |X(t)| \le u_T\right\} \to \exp\{-\exp(-x)\},$$

where $u = u_T$ is the solution to the equation

$$TH_2^{a/b} \frac{2}{\sqrt{2\pi}u} e^{-u^2/2} = e^{-x},$$

and where $H_2^R$ equals $\sqrt{1 + (1/R)}$ [cf. Konstant and Piterbarg (1993), page 87]. Solving this equation as, for example, in the proof of Theorem 1.5.3 in Leadbetter et al. (1983), gives the result. □

We will apply this theorem to the process

$$X_r(t) = Y^{(r)}(t)/\sigma_r$$

for $r = 3, 4$—which has period one and where $\sigma_r$ is given in Lemma 1—in combination with Proposition 5—to obtain the following result.

PROPOSITION 9. *Let $\mathcal{D}, j_n$ be as in Corollary 6. If $f_n^{(r)} := f_n(y, j_n)$ is the wavelet density estimator from (2.1) based on the quadratic ($r = 3$) or the cubic ($r = 4$) Battle–Lemarié scaling function $\phi_r$, then, for every $x \in \mathbb{R}$,*

$$\sup_{f \in \mathcal{D}} \left| \Pr_f \left\{ A_n \left( \sqrt{n 2^{-j_n}} \left\| \frac{f_n^{(r)} - Ef_n^{(r)}}{\sigma_r \sqrt{f}} \right\|_{[0,1]} - B_n^{(r)} \right) \le x \right\} - e^{-e^{-x}} \right| \to 0,$$

$$r = 3, 4,$$

*where $\sigma_r$ is given in Lemma 1, $A_n = A(j_n)$, $B_n^{(r)} = B^{(r)}(j_n)$ with*

(4.30)
$$A(l) = \sqrt{2(\log 2)l} \quad \text{and}$$

$$B^{(r)}(l) = A(l) - \frac{\log l + \log(\pi \log 2) - \log\sqrt{1 + D_r}}{A(l)}, \qquad r = 3, 4,$$

*$D_3 = \sum_k (a_k^{(3)} - a_{k-2}^{(3)})^2 / M$ with $M$ as in Lemma 1, and $D_4 = 9\sum_k (a_k^{(4)} - a_{k-2}^{(4)})^2 / C$ with $C$ as in (4.53).*

PROOF. We start with the case $r = 3$, the *quadratic* Battle–Lemarié wavelet. We first verify the hypotheses of Theorem 2 for the process $X_r$, with $\sigma_X(t) = \sigma(t)$ and $r_X(s,t) = r(s,t)$. We drop the sub- or superindex $r = 3$ from $X$, $Y$, $\sigma$. In this case Lemma 1 gives

$$\sigma^2 = \sigma^2(1/2) = \frac{1}{4}\sum_k (a_k + a_{k-1})^2 + \frac{M}{32}.$$

As in (4.23), we can write

$$Y(t) = \sum_{k \in \mathbb{Z}} \phi_3(t+k)g_k, \qquad t \in \mathbb{R},$$

where $g_k$ are i.i.d. standard normal and where

$$\phi_3(t+k) = \tfrac{1}{2}(a_k - 2a_{k-1} + a_{k-2})t^2 + (a_{k-1} - a_{k-2})t$$
$$+ \tfrac{1}{2}(a_{k-1} + a_{k-2}), \qquad 0 \le t \le 1,$$

by the computation in (4.49) below. It follows that $\phi_3$ is differentiable for all $t \in \mathbb{R}$ with the derivative

(4.31) $\quad \phi_3'(t+k) = (a_k - 2a_{k-1} + a_{k-2})t + (a_{k-1} - a_{k-2}), \qquad 0 \le t \le 1.$

In particular the process $Y'(t) = \sum_k \phi_3'(t+k)g_k$ is defined since the coefficients $a_k$ have exponential decay, and

$$E\left(\frac{Y(t+h) - Y(t) - hY'(t)}{h}\right)^2$$
$$= \sum_k \left(\frac{\phi_3(t+k+h) - \phi_3(t+k) - h\phi_3'(t+k)}{h}\right)^2 \to 0$$

as $h \to 0$ (for $0 \le t, t+h \le 1$, the quantity inside the parenthesis is dominated by $h|a_k - 2a_{k-1} + a_{k-2}|$ which is square summable). This shows that the process $Y$ (and hence also $X$) is differentiable in quadratic mean. Moreover, $Y$ (and $X$) has a.s. continuous sample paths [note that $Y(t) = t^2 G_1 + t G_2 + C$ for a constant $C$ and a bivariate normal variable $(G_1, G_2)$]. Hence, condition 1 in Theorem 2 is satisfied.

Condition 3 is also satisfied with $t_0 = 1/2$: By Lemma 1, $\sigma_X(1/2) = \sigma(1/2)/\sigma(1/2) = 1$, and this maximum is strict. Moreover, using (4.51) below and the comments after it, we have

$$\sigma^2(t) - \sigma^2(1/2) = \frac{M}{2}(t^4 - 2t^3 + t^2) - \frac{M}{32} = \frac{M}{2}\left[\left(t - \frac{1}{2}\right)^2 - \frac{1}{4}\right]^2 - \frac{M}{32}$$
$$= -\frac{M}{4}\left(t - \frac{1}{2}\right)^2 + \frac{M}{2}\left(t - \frac{1}{2}\right)^4,$$

which implies that $\sigma'(1/2) = 0$ and also gives

(4.32) $$\sigma''_X(1/2) = \sigma''(1/2)/\sigma = -M/(4\sigma^2).$$

Finally, using (4.31),

(4.33) $$E(X'(1/2))^2 = \sum_k (\phi'_3(1/2+k))^2/\sigma^2 = \sum_k (a_k - a_{k-2})^2/(4\sigma^2) > 0$$

(if the last sum were zero, the exponential decay of the $a_k$ would be contradicted) which completes the verification of condition (3) in Theorem 2.

We next consider condition 2. Recall (4.24) and (4.25); if this condition does not hold, then there exist $s$, $t$ in $[0,1)$ and $m \in \mathbb{N} \cup \{0\}$ with $s \neq t$ or with $m > 0$, or both, such that

$$\sum_k \phi_3(t+k)\phi_3(s+m+k) = \pm \sqrt{\sum_k \phi_3^2(t+k)} \sqrt{\sum_k \phi_3^2(s+k)}.$$

The right-hand side of this equation is different from zero for all $s$ and $t$ by (4.51) below. Consequently, this identity is satisfied if and only if there exist $\lambda \neq 0$ such that the vector with $k$th coordinate $\phi(t+k)$ equals $\lambda$ times the vector with $k$th coordinate $\phi(s+m+k)$. By (4.27) this is equivalent to the existence of $\lambda \neq 0$ satisfying the infinite system of equations

$$N_3(t)a_k + N_3(t+1)a_{k-1} + N_3(t+2)a_{k-2}$$
$$= \lambda N_3(s)a_{k+m} + \lambda N_3(s+1)a_{k+m-1} + \lambda N_3(s+2)a_{k+m-2},$$

$k \in \mathbb{Z}$, for some $s$, $t$ and $m$ satisfying the specified conditions. If we let $v^{(l)} \in \mathbb{R}^{\mathbb{Z}}$ be defined by the coordinates $v^{(l)}_k = a_{k-l}$, then we can write this system of equations as the vector equation

(4.34) $$N_3(t)v^{(0)} + N_3(t+1)v^{(1)} + N_3(t+2)v^{(2)}$$
$$= \lambda N_3(s)v^{(-m)} + \lambda N_3(s+1)v^{(-m+1)} + \lambda N_3(s+2)v^{(-m+2)}.$$

It is not difficult to see that if this equation has a solution with $m = 0$ and $s \neq t$, $s, t \in [0,1)$, or with $m > 0$ and $s, t \in [0,1)$, then a finite nonvoid set of vectors $v^{(l)}$ are linearly dependent. For instance, for $m = 0$, this equation becomes (using the explicit form of $N_3$ in Section 4.4.3)

$$(t^2 - \lambda s^2)v^{(0)} + 2[(3(1-\lambda)/4) - (t-1/2)^2 + \lambda(s-1/2)^2]v^{(1)}$$
$$+ ((1-t)^2 - \lambda(1-s)^2)v^{(2)} = 0,$$

and the coefficients are all zero if and only if $\lambda = 1$ and $s = t$ so that if this equation has a solution with $\lambda \neq 0$ and $s \neq t$, both in $[0,1)$, then the vectors $v^{(i)}$, $i = 0, 1, 2$, are linearly dependent. It is even easier to see that a similar conclusion holds for $m > 0$ (in this case, since $N_3(t+2) \neq 0$ for $t \in [0,1)$, if equation (4.34) has a solution for some $m > 0$, then the vectors $v^{(i)}, v^{(-m+i)}$, $i = 0, 1, 2$, are linearly dependent). However, suppose $\sum_{i=1}^{r} \lambda_i v^{(l_i)} = 0$ for some $r > 0$, some $l_i \in \mathbb{Z}$ and $\lambda_i \neq 0$. Since $\phi_3(\cdot - \ell) = \sum_k v_k^{(\ell)} N_3(\cdot - k)$, $\sum_{i=1}^{r} \lambda_i \phi_3(\cdot - l_i) = 0$ follows which is impossible unless $\lambda_i = 0$ for all $i$ because different integer translates of a father wavelet are orthogonal, hence linearly independent. This verifies condition 2 in Theorem 2.

Finally we check condition 4. Let us recall that, by the exponential decay of $\phi_3$, there exist $c_1, c_2 > 0$ such that $|\phi_3(x)| \leq c_1 \lambda^{-c_2|x|}$ for all $x$ so that, for $t - s \geq v > 0$,

$$\sigma^2 |r(s,t)| = \left| \sum_k \phi_3(s - k + (t-s))\phi_3(s-k) \right|$$

$$\leq \|\phi_3\|_\infty \left( \sum_{k: |k-s| \leq (t-s)/2} |\phi_3(s - k + (t-s))| \right.$$

$$\left. + \sum_{k: |k-s| > (t-s)/2} |\phi_3(s-k)| \right)$$

$$\leq 2c_1 \|\phi_3\|_\infty \lambda^{-c_2 v/2},$$

proving condition 4. Hence, the process $X(t) = Y(t)/\sigma$ satisfies the hypotheses of Theorem 2 with $T = 2^{j_n}$ and constants given in (4.32) and (4.33). Combining the resulting asymptotic distribution for $\sup_{0 \leq t \leq 2^{j_n}} |Y(t)|$ with Proposition 5 (which applies, as shown in Section 4.3.1), we obtain Proposition 9 for the quadratic Battle–Lemarié wavelet.

Now we consider the cubic Battle–Lemarié wavelet density estimator (case $r = 4$). Theorem 2 cannot be applied directly to $X_4(t) = \sum_k \phi_4(t-k)g_k/\sigma_4$ because the variance of this cyclostationary Gaussian process of period 1 has its maxima at $t = k$ which contradicts condition 3 ($t_0 = 0, 1$ are not unique and are not interior to $[0,1]$). This is not a serious difficulty and can be easily dealt with. With $Y_n(t) := Y_n^{(4)}(t) = \sum_k \phi_4(2^{j_n} t + k)g_k$, we will obtain a limit theorem for

$$\sup_{t \in [0,1]} |Y_n(t)|/\sigma = \sup_{t \in [0, 2^{j_n}]} \left| \sum_k \phi_4(t+k)g_k \right| \Big/ \sigma := \sup_{t \in [0, 2^{j_n}]} |Y(t)|/\sigma$$

via limit theorems for

$$Z_n^- = \sup_{\delta/2^{j_n} \leq t \leq 1 - \delta/2^{j_n}} |Y_n(t)|/\sigma \quad \text{and} \quad Z_n^+ = \sup_{-\delta/2^{j_n} \leq t \leq 1 + \delta/2^{j_n}} |Y_n(t)|/\sigma,$$

where

$$\sigma^2 = \sigma_4^2 = \frac{1}{36} \sum_k (a_{k-1} + 4a_k + a_{k+1})^2$$

(cf. Lemma 1) and where $0 < \delta < 1$ (e.g., $\delta = 1/2$). Set $Z(t) = Y(t+\delta)/\sigma$, $t \geq 0$. Then,

$$Z_n^- = \sup_{0 \leq t \leq 2^{j_n} - 2\delta} |Z(t)|$$

and $Z(t)$ is still a cyclostationary centered Gaussian process since $r_Z(t, t+v) = r_Y(t+\delta, t+\delta+v)$ which is periodical with period 1 for each $v$. But now $\sigma_Z^2(t)$ attains a unique and strict maximum on $[0,1]$ at the point $t_0 = 1 - \delta$. It is easy to see, proceeding in analogy with the quadratic spline case and using computations from Section 4.4.3, that conditions 1–4 in Theorem 2 are satisfied with

$$\sigma_Z''(1-\delta) = -C/(36\sigma^2) \quad \text{and} \quad E(Z'(1-\delta))^2 = \sum_k (a_k - a_{k-2})^2/(4\sigma^2),$$

where $C > 0$ is defined in (4.53). Therefore, if we set

$$A_n^- = \sqrt{2\log(2^{j_n} - 2\delta)},$$

$$B_n^- = A_n^- - \frac{\log\log(2^{j_n} - 2\delta) + \log\pi - \log\sqrt{1+D_4}}{A_n^-},$$

Theorem 2 proves that the random variables $A_n^-(Z_n^- - B_n^-)$ converge in law to the Gumbel distribution. Now, let $A_n$ and $B_n = B_n^{(4)}$ be the constants in (4.30). We have

$$A_n(Z_n^- - B_n) = A_n^-(Z_n^- - B_n^-)\frac{A_n}{A_n^-} + A_n(B_n^- - B_n),$$

and, as is easy to see,

$$\frac{A_n^-}{A_n} \to 1 \quad \text{and} \quad A_n(B_n^- - B_n) \to 0 \quad \text{as } n \to \infty.$$

We thus conclude that the sequence $\{A_n(Z_n^- - B_n)\}$ is weak convergence equivalent to the sequence $\{A_n^-(Z_n^- - B_n^-)\}$, and therefore, for all $x \in \mathbb{R}$,

$$\Pr\{A_n(Z_n^- - B_n) \leq x\} \to e^{-e^{-x}}.$$

The same argument gives (with $Z(t) = Y(t-\delta)/\sigma$)

$$\Pr\{A_n(Z_n^+ - B_n) \leq x\} \to e^{-e^{-x}};$$

hence, since $Z_n^- \leq \sup_{t \in [0,1]} |Y_n(t)|/\sigma \leq Z_n^+$, we have proved that for all $x$,

$$\Pr\left\{A_n\left(\sup_{t \in [0,1]} |Y_n(t)|/\sigma - B_n\right) \leq x\right\} \to e^{-e^{-x}}.$$

Combining this limit with Proposition 5 which, as argued in Section 4.3.1, applies to the projection kernels of the cubic spline wavelet, we obtain Proposition 9 for $r = 4$. □

REMARK 3. Regarding higher order Battle–Lemarié wavelets ($r > 4$), notice that—the scaling function $\phi_r$ being a piecewise polynomial function with smooth (and nonconstant) weldings—the absolute maximum of the variance function (4.25) on $[0,1]$ is attained at a finite number of points, perhaps even at a single point (as in the cases $r \leq 4$). In this case results from the literature mentioned before Theorem 2 can be applied, and proofs can be given along the lines of the proof of Proposition 9. After this paper was written it was shown that the variance function (4.25) has unique maxima in the case $r \leq 9$ [see Giné and Madych (2009), where one can also find a more general conjecture related to these questions].

### 4.4. Remaining proofs.

4.4.1. *Proof of Theorem 1.* Throughout this proof we shall often write, in slight abuse of notation, $f_n(l)$ for $f_n(\cdot, l)$. Note also that all densities in $\mathcal{P}$ are bounded by a fixed constant depending only on $\eta, b$, and we shall denote this constant by $U$. For every $f \in \mathcal{P}$ there exists a unique $t := t(f)$ such that $f$ satisfies Condition 3 for this $t$. Define

$$B(j,t) = b_2 2^{-jt}, \qquad \sigma(n,l) = \sqrt{\frac{2^l l}{n}},$$

(4.35)
$$j_n^*(t) = \min\left\{ j \in \mathcal{J} : B(j,t) \leq \frac{M\sigma(n_2, j)}{4} \right\}.$$

It is easy to see that $j_n^*(t)$ satisfies

$$(4.36) \qquad 2^{j_n^*(t)} \simeq \left(\frac{n_2}{\log n_2}\right)^{1/(2t+1)},$$

so it is a "rate optimal" resolution level for estimating $f \in \mathcal{C}^t(\mathbb{R})$. We wish to show that $\hat{j}_n$ defined in (3.1) is a 'good estimate' for $j_n^*(t)$, which is achieved by the following lemma. Recall that $\hat{j}_n$ is defined so far only up to the choice of the constant $M'$ (cf. Remark 2), which can be retrieved from the proof of the following lemma (and Proposition 10).

LEMMA 2. (a) *We can choose a finite positive constant $M'$ depending only on $K$ such that if $j \in \mathcal{J}_n$ satisfies $j > j_n^*(t)$; then for every $n \in \mathbb{N}$,*

$$\Pr_f(\hat{j}_n = j) \leq c 2^{-j/c},$$

*where the constant $0 < c < \infty$ depends only on $M'$ and $K$.*

(b) *There exists a positive integer $m$ depending only on $b_1, b_2, \eta$, and a constant $c'$ depending only on $M'$ from part* (a) *and on $K$ such that for every $j \in \mathcal{J}_n$ satisfying $j < j_n^*(t) - m$ and every $n \in \mathbb{N}$ large enough such that $j_n^*(t) \geq 2$, we have*

$$\Pr_f(\hat{j}_n = j) \leq c' 2^{-j/c'}.$$

PROOF. Since this lemma only involves the sample $S_2$, we set $n = n_2$ for notational simplicity. To prove the first claim, fix $j > j_n^*(t)$ and let $j^- = j - 1$. Then one has

$$(4.37) \quad \Pr_f(\hat{j}_n = j) \leq \sum_{l \in \mathcal{J}: l \geq j} \Pr_f(\| f_n(j^-) - f_n(l) \|_\infty > M\sigma(n, l)).$$

We first observe that by Condition 3,

$$\| f_n(j^-) - f_n(l) \|_\infty$$
$$\leq \| f_n(j^-) - f_n(l) - Ef_n(j^-) + Ef_n(l) \|_\infty + B(j^-, t) + B(l, t)$$

and that

$$B(j^-, t) + B(l, t) \leq 2B(j_n^*(t), t) \leq (M/2)\sigma(n, j_n^*(t)) \leq (M/2)\sigma(n, l)$$

by definition of $j_n^*(t)$ and since $l > j^- \geq j_n^*(t)$. Consequently, the $l$th probability in the last sum is bounded by

$$\Pr_f(\| f_n(j^-) - f_n(l) - Ef_n(j^-) + Ef_n(l) \|_\infty > (M - (M/2))\sigma(n, l))$$
$$\leq \Pr_f(\| f_n(j^-) - Ef_n(j^-) \|_\infty > (M/4)\sigma(n, l))$$
$$+ \Pr_f(\| f_n(l) - Ef_n(l) \|_\infty > (M/4)\sigma(n, l)) \leq d2^{-l/d},$$

where $0 < d < \infty$ depends only on $M', K$, and where we have used the fact that we can choose $M'$ large enough but finite depending only on $K$ so that Proposition 10 below applies. Feeding this bound into (4.37) and summing the series proves the first claim of the lemma.

To prove the second claim, fix $j < j_n^*(t) - m$ where $m$ is some positive integer, and observe that

$$(4.38) \quad \Pr_f(\hat{j}_n = j) \leq \Pr_f(\| f_n(j) - f_n(j_n^*(t)) \|_\infty \leq M\sigma(n, j_n^*(t))).$$

Now using Condition 3 and the triangle inequality we deduce

$$\| f_n(j) - f_n(j_n^*(t)) \|_\infty \geq (b_1/b_2)B(j, t) - B(j_n^*(t), t)$$
$$- \| f_n(j) - Ef_n(j) - f_n(j_n^*(t)) + Ef_n(j_n^*(t)) \|_\infty$$

so that the probability in (4.38) is bounded by

$$\Pr_f(\| f_n(j) - Ef_n(j) - f_n(j_n^*(t)) + Ef_n(j_n^*(t)) \|_\infty$$
$$\geq (b_1/b_2)B(j, t) - B(j_n^*(t), t) - M\sigma(j_n^*(t), n)).$$

By definition of $j_n^*(t)$ and $B(j, t)$, we have

$$\frac{b_1}{b_2}B(j, t) - B(j_n^*(t), t) = \frac{b_1}{b_2}2^{t(j_n^*(t) - j)}B(j_n^*(t), t) - B(j_n^*(t), t)$$
$$> \left(\frac{b_1}{b_2}2^{t(m-1)} - 2^{-t}\right)B(j_n^*(t) - 1, t)$$
$$\geq \left(\frac{b_1}{b_2}2^{\eta(m-1)} - 2^{-\eta}\right)B\left(j_n^*(t) - 1, t\right),$$

and, using also $\sqrt{(j_n^*(t)-1)/j_n^*(t)} \geq 1/\sqrt{2}$ in view of $j_n^*(t) \geq 2$,
$$B(j_n^*(t)-1,t) \geq (M/4)\sigma(n, j_n^*(t)-1) = (M/4)2^{-1}\sigma(n, j_n^*(t))$$
so that the last probability is bounded by
$$\Pr_f\bigg(\|f_n(j) - Ef_n(j) - f_n(j_n^*(t)) + Ef_n(j_n^*(t))\|_\infty$$
$$\geq \bigg[\frac{M}{8}\bigg(\frac{b_1}{b_2}2^{\eta(m-1)} - 2^\eta\bigg) - M\bigg]\sigma(n, j_n^*(t))\bigg)$$
$$\leq \Pr_f\bigg(\|f_n(j) - Ef_n(j)\|_\infty \geq 2^{-1}\bigg[\frac{M}{8}\bigg(\frac{b_1}{b_2}2^{\eta(m-1)} - 2^\eta\bigg) - M\bigg]\sigma(n, j)\bigg)$$
$$+ \Pr_f\bigg(\|f_n(j_n^*(t)) - Ef_n(j_n^*(t))\|_\infty$$
$$\geq 2^{-1}\bigg[\frac{M}{8}\bigg(\frac{b_1}{b_2}2^{\eta(m-1)} - 2^\eta\bigg) - M\bigg]\sigma(n, j_n^*(t))\bigg).$$

We can now choose $m > 2$ sufficiently large but finite and only depending on $b_1, b_2, \eta$ so that, using Proposition 10 below, the last two probabilities can be made less than $c' 2^{-j/c'}$ for some constant $c'$ that depends only on $M', K$. $\square$

For the rest of the proof we assume that $M'$ and $n$ have been chosen large enough so that Lemma 2 holds. As a first consequence we obtain

$$\sup_{f\in\mathcal{P}} \Pr_f\{2^{\hat{j}_n} > 2^{j_n^*(t)}\} = \sup_{f\in\mathcal{P}} \Pr_f\{\hat{j}_n > j_n^*(t)\}$$

(4.39)
$$\leq \sup_{f\in\mathcal{P}} \sum_{j_n^*(t) < j \leq j_{\max}} \Pr_f(\hat{j}_n = j)$$
$$\leq c \sum_{j_n^*(t) < j \leq j_{\max}} 2^{-j/c} \leq c'' 2^{-j_n^*(t)/c''} \to 0$$

as $n \to \infty$ which already proves (3.8) by using (4.36) and the definition of $\hat{\sigma}_n$. [Note that $\hat{j}_n \in \mathcal{J}_n$ for all $n$ implies that we have to prove (3.8) only for large $n$.] Moreover,

$$\sup_{f\in\mathcal{P}} \Pr_f(2^{-\hat{j}_n} > 2^m 2^{-j_n^*(t)}) \leq \sup_{f\in\mathcal{P}} \Pr_f(\hat{j}_n < j_n^*(t) - m)$$

(4.40)
$$= \sup_{f\in\mathcal{P}} \sum_{j_{\min} \leq j < j_n^*(t)-m} \Pr_f(\hat{j}_n = j)$$
$$\leq c' \sum_{j_{\min} \leq j < j_n^*(t)-m} 2^{-j/c'} = c''' 2^{-j_{\min}/c'''} \to 0$$

as $n \to \infty$. Combining the above arguments we also have

(4.41) $$\sup_{f \in \mathcal{P}} \Pr_f\{\hat{j}_n \notin [j_n^*(t) - m, j_n^*(t)]\} \to 0$$

as $n \to \infty$ which we shall use in the following lemma.

LEMMA 3. *We have (E denoting expectation w.r.t. $S_1$)*

$$\left| \hat{\sigma}_n^{-1} \left\| \frac{f_{n_1}(\hat{j}_n + u_n) - f}{c(K)\sqrt{f_{n_1}(\hat{j}_n + u_n)}} \right\|_{[0,1]} - \hat{\sigma}_n^{-1} \left\| \frac{f_{n_1}(\hat{j}_n + u_n) - Ef_{n_1}(\hat{j}_n + u_n)}{c(K)\sqrt{f}} \right\|_{[0,1]} \right|$$
$$= o_P(1/\sqrt{\log n})$$

*uniformly in $\mathcal{P}$.*

PROOF. We will use $n \simeq n_1 \simeq n_2$ without mentioning in this proof. A few preliminary observations are necessary: We first establish that

(4.42) $$\hat{\sigma}_n^{-1} \| f_{n_1}(\hat{j}_n + u_n) - f \|_\infty = O_P(\sqrt{\log n})$$

uniformly in $\mathcal{P}$; note that by definition of $j_n^*(t)$ we have

$$B(j + u_n, t) \le (M/4)\sigma(n_2, j + u_n) \le C''\sigma(n_1, j + u_n)$$

for $j \ge j_n^*(t) - m$ and for $n$ large enough s.t. $u_n > m$ and where $C''$ depends only on $K$, $n_1/n_2$ and the envelope $U$ of $\mathcal{P}$. Hence, using (4.41) and Proposition 10,

$$\Pr_f\{\hat{\sigma}_n^{-1}\|f_{n_1}(\hat{j}_n + u_n) - f\|_\infty \ge C\sqrt{\log n}\}$$
$$\le \sum_{j \in [j_n^*(t)-m, j_n^*(t)]} \Pr_f\{\sigma^{-1}(n_1, j+u_n)\|f_{n_1}(j+u_n) - f\|_\infty \ge CC'\} + o(1)$$
$$\le \sum_{j \in [j_n^*(t)-m, j_n^*(t)]} \Pr_f\{\|f_{n_1}(j+u_n) - Ef_{n_1}(j+u_n)\|_\infty$$
$$\ge CC'\sigma(n_1, j+u_n) - B(j+u_n, t)\} + o(1)$$
$$\le \sum_{j \in [j_n^*(t)-m, j_n^*(t)]} \Pr_f\{\|f_{n_1}(j+u_n) - Ef_{n_1}(j+u_n)\|_\infty$$
$$\ge (CC' - C'')\sigma(n_1, j+u_n)\} + o(1)$$
$$= o(1)$$

for $C$ and $n$ large enough depending only on $\mathcal{P}$ (so that this bound is uniform in $\mathcal{P}$).

Second, (4.36), (4.39), (4.42) and Condition 4 imply that

$$\|f_{n_1}(\hat{j}_n + u_n) - f\|_\infty = O_P\left(\sqrt{\frac{2^{j_n^*(t)}\log n}{n}}2^{u_n/2}\right) = o_P(1/\log^2 n)$$

uniformly in $\mathcal{P}$, and since $f \geq \delta$ on $[0,1]$, we also have

(4.43) $$\sup_{y\in[0,1]} |f_{n_1}^{-1}(y, \hat{j}_n + u_n)| = O_P(1),$$

uniformly in $\mathcal{P}$. Consequently, since the root-transformation is Lipschitz on intervals bounded away from zero, we obtain, uniformly in $\mathcal{P}$, $\|f_{n_1}^{-1/2}(\hat{j}_n + u_n)\|_{[0,1]} = O(1)$ and

(4.44) $$\sup_{y\in[0,1]} |\sqrt{f_{n_1}(y, \hat{j}_n + u_n)} - \sqrt{f(y)}| = o_P(1/\log^2 n).$$

We now prove the lemma. First, using the above facts,

$$\left|\hat{\sigma}_n^{-1}\left\|\frac{f_{n_1}(\hat{j}_n + u_n) - f}{c(K)\sqrt{f_{n_1}(\hat{j}_n + u_n)}}\right\|_{[0,1]} - \hat{\sigma}_n^{-1}\left\|\frac{f_{n_1}(\hat{j}_n + u_n) - f}{c(K)\sqrt{f}}\right\|_{[0,1]}\right|$$

$$\leq \hat{\sigma}_n^{-1}\left\|\frac{(f_{n_1}(\hat{j}_n + u_n) - f)(\sqrt{f_{n_1}(\hat{j}_n + u_n)} - \sqrt{f})}{c(K)\sqrt{f}\sqrt{f_{n_1}(\hat{j}_n + u_n)}}\right\|_{[0,1]} = o_P(1/\log n).$$

Second, using Condition 3, (4.40) and Condition 4, we have for some constant $d$ that depends only on $b_2, c(K), \delta$ that

$$\left|\hat{\sigma}_n^{-1}\left\|\frac{f_{n_1}(\hat{j}_n + u_n) - f}{c(K)\sqrt{f}}\right\|_{[0,1]} - \hat{\sigma}_n^{-1}\left\|\frac{f_{n_1}(\hat{j}_n + u_n) - Ef_{n_1}(\hat{j}_n + u_n)}{c(K)\sqrt{f}}\right\|_{[0,1]}\right|$$

$$\leq \hat{\sigma}_n^{-1}\left\|\frac{Ef_{n_1}(\hat{j}_n + u_n) - f}{c(K)\sqrt{f}}\right\|_{[0,1]}$$

$$\leq d\hat{\sigma}_n^{-1} 2^{-t(\hat{j}_n + u_n)}$$

$$= O_P(\sqrt{n} 2^{-j_n^*(t)(t+(1/2))} 2^{-u_n(t+(1/2))})$$

$$= o_P(1/\sqrt{\log n})$$

since $f$ is bounded from below and since $2^{-j_n^*(t)(t+1/2)} \simeq (n/\log n)^{-1/2}$ by (4.36). □

By the above lemma, and since $\hat{A}_n = A(\hat{j}_n + u_n) \simeq \sqrt{\log n}$, to complete the proof of the theorem, it suffices to prove the required limit for

$$\Pr_f\left\{A(\hat{j}_n + u)\left(\hat{\sigma}_n^{-1}\left\|\frac{f_{n_1}(\hat{j}_n + u_n) - Ef_{n_1}(\hat{j}_n + u_n)}{c(K)\sqrt{f}}\right\|_{[0,1]} - B(\hat{j}_n + u_n)\right) \leq x\right\}$$

$$= \Pr{}_f(W_n(\hat{j}_n + u_n, x)),$$

where we introduce the events

$$W_n(j, x) = \left\{ \sqrt{\frac{n_1}{2^j}} \left\| \frac{f_{n_1}(j) - Ef_{n_1}(j)}{c(K)\sqrt{f}} \right\|_{[0,1]} \le B(j) + \frac{x}{A(j)} \right\}$$

for $j \in \mathbb{N}$. Recall that Condition 2 implies, for every $j_n \in \mathcal{J}$ that

$$\Pr{}_f(W_n(j_n, x)) \to \zeta(x) := \exp\{-\exp\{-x\}\}$$

as $n \to \infty$, uniformly in $\mathcal{P}$, and note that $j + u_n \in \mathcal{J}_n$ for $j \in [j_n^*(t) - m, j_n^*(t)]$ and $n$ large enough by (4.36). Using (4.41) and independence of $\hat{j}_n$ and $f_{n_1}(j)$, we have

$$\Pr{}_f(W_n(\hat{j}_n + u_n, x))$$
$$= \sum_{j \in \mathcal{J}} \Pr{}_f(W_n(j + u_n, x) \cap \{\hat{j}_n = j\})$$
$$= \sum_{j_n^*(t) - m \le j \le j_n^*(t)} \Pr{}_f(W_n(j + u_n, x) \cap \{\hat{j}_n = j\}) + o(1)$$
$$= \sum_{j_n^*(t) - m \le j \le j_n^*(t)} \Pr{}_f(W_n(j + u_n, x)) \Pr{}_f(\{\hat{j}_n = j\}) + o(1)$$
$$= \zeta(x) \sum_{j_n^*(t) - m \le j \le j_n^*(t)} \Pr{}_f(\hat{j}_n = j)$$
$$+ \sum_{j_n^*(t) - m \le j \le j_n^*(t)} (\Pr{}_f(W_n(j, x)) - \zeta(x)) \Pr{}_f(\hat{j} = j) + o(1).$$

But the limit of the last expression is $\zeta(x)$, completing the proof of Theorem 1. The first quantity converges to $\zeta(x)$ as $n \to \infty$ since

$$\sum_{j_n^*(t) - m \le j \le j_n^*(t)} \Pr{}_f(\hat{j}_n = j) = 1 - \sum_{j \notin [j_n^*(t) - m, j_n^*(t)]} \Pr{}_f(\hat{j}_n = j) \to 1$$

uniformly in $f \in \mathcal{P}$ by (4.41). The second quantity converges to zero since

$$\max_{j_n^*(t) - m \le j \le j_n^*(t)} |\Pr{}_f(W_n(j, x)) - \zeta(x)| \to 0,$$

in view of Condition 2 and since $m$ is finite (depending only on $\mathcal{P}$).

4.4.2. *Proofs and complementary results for Section 3.5.*

PROOF OF PROPOSITION 4. We use the fact that the wavelet basis $\phi, \psi$ characterizes the space $\mathcal{C}^t(\mathbb{R})$ for $t < r - 1$, cf. Definition 1. Since

$$|\beta_{lk}(g)| = \left| 2^{l/2} \int_\mathbb{R} \psi(2^l x - k) g(x) \, dx \right|$$

$$\leq 2^{-l/2}\|\psi\|_1\|g\|_\infty$$

for every $l,k$ and every bounded function $g$, we have for $g = K_j(f) - f$, whose wavelet coefficients are zero for $l < j$, that

$$(4.45) \qquad \|K_j(f) - f\|_\infty \geq \|\psi\|_1^{-1} \sup_{l \geq j, k \in \mathbb{Z}} |2^{l/2}\beta_{lk}(f)|.$$

Take

$$E_m(k) = \{f \in \mathcal{C}^t(\mathbb{R}) : |\beta_{lk}(f)| \geq 2^{-l(t+1/2)}2^{-m} \text{ for every } l \in \mathbb{N}\}$$

and set

$$A_m(k) = \{h \in \mathcal{C}^t(\mathbb{R}) : \|h - f\|_{t,\infty} < 2^{-m-1} \text{ for some } f \in E_m(k)\}$$

so that (recalling Definition 1), for every $h \in A_m(k)$ and every $l$, $|\beta_{lk}(h)| \geq 2^{-l(t+1/2)}2^{-m-1}$. Consequently, using (4.45), we have

$$(4.46) \qquad \|K_j(h) - h\|_\infty \geq \|\psi\|_1^{-1} 2^{-m-1} 2^{-jt}$$

for every nonnegative integer $j$ and every $h \in A_m(k)$. Define now

$$A = \bigcup_{m \geq 0, k \in \mathbb{Z}} A_m(k),$$

all of whose elements satisfy the lower bound (4.46) for some $m$, and therefore $A \subset \mathcal{N}_t^c$. $A$ is clearly open and it is also dense in $\mathcal{C}^t(\mathbb{R})$: Let $g \in \mathcal{C}^t(\mathbb{R})$ be arbitrary, and define the function $g_m$ by its wavelet coefficients $\alpha_k(g_m) = \alpha_k(g)$ for all $k$ and $\beta_{lk}(g_m)$ equal to $\beta_{lk}(g)$ when $|\beta_{lk}(g)| > 2^{-l(t+1/2)}2^{-m}$ and equal to $2^{-l(t+1/2)}2^{-m}$, otherwise. Clearly $g_m \in A$ for every $m$, and for $\varepsilon > 0$ arbitrary we can choose $m$ large enough such that

$$\|g - g_m\|_{t,\infty} = \sup_{l,k} 2^{l(t+1/2)}|\beta_{lk}(g) - 2^{-l(t+1/2)}2^{-m}|1_{|\beta_{lk}(g)| \leq 2^{-l(t+1/2)}2^{-m}}$$

$$\leq 2^{-m+1} < \varepsilon.$$

This proves that $\mathcal{N}_t$ is contained in the complement of an open and dense set, hence itself must be nowhere dense. $\square$

We now construct some explicit densities in the class $\mathcal{P}$, further illustrating the genericity of Condition 3. For the case of compactly supported wavelets, let $f$ be any density in the class $\mathcal{D}(t, D, \delta, F)$. The wavelet series of $f$ is

$$f = K_0(f) + \sum_{l=0}^{\infty} \sum_k \beta_{lk}(f)\psi_{lk}.$$

If $\psi$ is supported in $[0,a]$ (if necessary after a translation), pick a $k_0 \in \mathbb{Z}$ and $l_0$ large so that $[k_0/2^l, (k_0+a)/2^l] \in [0,1]$ for every $l \geq l_0$, and define $g(x)$ to have exactly the same wavelet series as $f$, but, for $l \geq l_0$, $\beta_{lk_0}(g) = \beta_{lk_0}(f)$ if $|\beta_{lk_0}(f)| \geq 2^{-l(t+1/2)}$ and $\beta_{lk_0}(g) = 2^{-l(t+1/2)}$, otherwise. By choice of $k_0$ this modification takes place only in $[k_0/2^l, (k_0+a)/2^l] \in [0,1]$ where $f$ is larger than or equal to $\delta$. Moreover,

$$\|f - g\|_\infty = \left\| \sum_{l \geq l_0} (\beta_{lk_0}(f) - 2^{-l(t+1/2)}) 1\{|\beta_{lk_0}(f)| \leq 2^{-l(t+1/2)}\} \psi_{lk_0} \right\|_\infty$$
$$\leq 2\|\psi\|_\infty \sum_{l \geq l_0} 2^{-lt}$$

so that by choosing $l_0$ large enough this quantity can be made as small as desired; in particular $g(x) \geq \delta(1 - \varepsilon)$ for $x \in [F_1, F_2]$ and every $\varepsilon > 0$. Furthermore, since $f$ integrates to one, and since $\psi$ integrates to zero, $g$ is a density, and its wavelet coefficients at $k_0$ satisfy the lower bound $|\beta_{lk_0}(g)| \geq 2^{-l(t+1/2)}(1-\varepsilon)$ for $l \geq l_0$. Using (4.45) this verifies the lower bound in Condition 3.

If $\psi$ does not have compact support, as is the case for Battle–Lemarié wavelets, the above modification might lead to functions that are negative somewhere, and hence not densities. To remedy this, start with any function $\bar{f}$ that is in $\mathcal{C}^t(\mathbb{R})$ and bounded from below by $\delta$ for *every* $x \in \mathbb{R}$. Then the above modification at $k_0$, $l \geq l_0$ (and the proof) gives us a function $f$ which is bounded from below by $\delta/2$ on all of $\mathbb{R}$ and satisfies $|\beta_{lk_0}(f)| \geq 2^{-l(t+1/2)}$ for every $l \geq l_0$. Multiply $f$ by a positive, bounded, integrable, infinitely-differentiable function $h$ (possibly compactly-supported) that is equal to one on $[-1, 1]$, and then divide $fh$ by $\|fh\|_1$, giving a (possibly compactly supported) density on $\mathbb{R}$ which is again contained in $\mathcal{C}^t(\mathbb{R})$ since this space is a multiplication algebra. The wavelet coefficients at $k_0$ are

$$|\beta_{lk_0}(fh)| = \left| \int_\mathbb{R} f(x)h(x)\psi_{lk_0}(x)\,dx \right|$$
$$= \left| \int_{-1}^1 f(x)\psi_{lk_0}(x)\,dx + \int_{[-1,1]^c} f(x)h(x)\psi_{lk_0}(x)\,dx \right|$$
$$\geq |\beta_{lk_0}(f)| - 2^{l/2}(\|f\|_\infty(1 + \|h\|_\infty)) \int_{[-1,1]^c} |\psi(2^l x - k_0)|\,dx,$$

and the quantity we subtract is, by exponential decay of $\psi$, less than or equal to a constant times $2^{-l/2} \lambda^{c|2^l - k_0|}$ for some $0 < \lambda < 1$. For $l \geq l_0$ large enough this quantity can be made smaller than any power of $2^l$, hence the same lower bound for the wavelet coefficients at $k_0$ holds for the density $fh$. Again, using (4.45) this verifies the lower bound in Condition 3.

4.4.3. *Proof of Lemma 1.* *The piecewise linear Battle–Lemarié wavelet.* In this case

$$N_2(x) = I_{[0,1)} * I_{[0,1)}(x) = \begin{cases} x, & \text{for } 0 \leq x \leq 1, \\ 2-x, & \text{for } 1 \leq x \leq 2, \end{cases}$$

which yields (writing $a_k$ for $a_k^{(2)}$)

(4.47)
$$\phi_2(k+t) = a_k N_2(t) + a_{k-1} N_2(t+1)$$
$$= a_k t + a_{k-1}(1-t), \quad k \in \mathbb{Z}, 0 \leq t \leq 1.$$

Then,

$$\sum \phi_2^2(k) = \sum a_k^2$$

and, using $2\sum a_{k-1}(a_k - a_{k-1}) = -\sum(a_k - a_{k-1})^2$,

(4.48)
$$\sum \phi_2^2(k+t) = \sum(a_k - a_{k-1})^2 t^2 + 2\sum a_{k-1}(a_k - a_{k-1})t + \sum \phi_2^2(k)$$
$$= \sum(a_k - a_{k-1})^2 t(t-1) + \sum \phi_2^2(k).$$

Now, $\sum(a_k - a_{k-1})^2 > 0$; otherwise all the $a_k$ would be identical which is impossible because the $a_k$ decay exponentially. Therefore, $\sum_k \phi_2(k+t)$ has a unique maximum on $[0,1)$, at $t=0$, that is, the variance function $\sum \phi^2(t-k)$ has only isolated maxima which are at the points $t = l \in \mathbb{Z}$.

*The quadratic Battle–Lemarié wavelet.* In this case,

$$N_3(x) = \begin{cases} x^2/2, & \text{for } 0 \leq x < 1, \\ 3/4 - (x - 3/2)^2, & \text{for } 1 \leq x < 2, \\ (3-x)^2/2, & \text{for } 2 \leq x \leq 3. \end{cases}$$

We then have (still omitting the superindex 3 from $a_k^{(3)}$)

$$\phi_3(k) = a_{k-1} N_3(1) + a_{k-2} N_3(2) = (a_{k-1} + a_{k-2})/2 \quad \text{and}$$

$$\sum \phi_3^2(k) = \frac{1}{4}\sum(a_k + a_{k-1})^2,$$

and, for $0 \leq t \leq 1$,

(4.49)
$$\phi_3(k+t) = a_k N_3(t) + a_{k-1} N_3(t+1) + a_{k-2} N_3(t+2)$$
$$= a_k t^2/2 + a_{k-1}[3/4 - (t-1/2)^2] + a_{k-2}(1-t)^2/2$$
$$= \left[\frac{a_k - a_{k-1}}{2} - \frac{a_{k-1} - a_{k-2}}{2}\right]t^2 + [a_{k-1} - a_{k-2}]t$$
$$+ \frac{1}{2}[a_{k-1} + a_{k-2}].$$

Hence,

$$\sum \phi_3^2(k+t) = \sum \left[\frac{a_k - a_{k-1}}{2} - \frac{a_{k-1} - a_{k-2}}{2}\right]^2 t^4$$
$$+ \sum [(a_k - a_{k-1}) - (a_{k-1} - a_{k-2})](a_{k-1} - a_{k-2})t^3$$
$$+ \sum \left[(a_{k-1} - a_{k-2})^2 \right.$$
$$\left. + \frac{1}{2}[(a_k - a_{k-1}) - (a_{k-1} - a_{k-2})](a_{k-1} + a_{k-2})\right] t^2$$
$$+ \sum (a_{k-1} - a_{k-2})(a_{k-1} + a_{k-2})t + \sum \phi_3^2(k).$$

But

$$\sum (a_{k-1} - a_{k-2})(a_{k-1} + a_{k-2}) = 0,$$
$$\sum (a_k - a_{k-1})a_{k-1} = -\frac{1}{2}\sum(a_k - a_{k-1})^2$$

and

$$\sum \frac{1}{2}(a_k - a_{k-1})(a_{k-1} + a_{k-2})$$
$$= \frac{1}{2}\sum(a_k - a_{k-1})(a_{k-1} - a_{k-2}) + \sum(a_k - a_{k-1})a_{k-2}$$
$$= \frac{1}{2}\sum(a_k - a_{k-1})(a_{k-1} - a_{k-2})$$
$$- \sum(a_k - a_{k-1})(a_{k-1} - a_{k-2}) + \sum(a_k - a_{k-1})a_{k-1}$$
$$= -\frac{1}{2}\sum(a_k - a_{k-1})(a_{k-1} - a_{k-2}) - \frac{1}{2}\sum(a_k - a_{k-1})^2.$$

So, setting

(4.50) $$M = \sum(a_k - a_{k-1})^2 - \sum(a_k - a_{k-1})(a_{k-1} - a_{k-2})$$

we obtain

(4.51)
$$\sum \phi_3^2(k+t) = \frac{1}{2}Mt^4 - Mt^3 + \frac{1}{2}Mt^2 + \sum \phi_3^2(k)$$
$$= \frac{1}{2}Mt^2(t-1)^2 + \sum \phi_3^2(k), \qquad 0 \le t \le 1.$$

We can write $M$ as $M = \|u\|\|v\| - \langle u, v\rangle$ where $\|u\| = \|v\|$, and so $M = 0$ iff $u = v$, but this would mean

$$\cdots = a_1 - a_0 = a_0 - a_{-1} = a_{-1} - a_{-2} = \cdots,$$

that is, that the points $(k, a_k)$ lie on a straight line which contradicts the exponential decay of the $a_k$'s. So $M$ is strictly positive. Therefore, $\sum \phi_3^2(t-k)$, $t \in [0,1]$, has a unique maximum at $t = 1/2$, $\sum \phi_3^2(k+1/2) = \frac{1}{32}M + \sum \phi_3^2(k)$. That is, the variance function has only isolated maxima which are at the points $k + 1/2$, $k \in \mathbb{Z}$.

*The cubic Battle–Lemarié wavelet.* We have [e.g., Schumaker (1993), page 136]

$$N_4(x) = \begin{cases} x^3/6, & \text{for } 0 \leq x \leq 1, \\ \frac{2}{3} - \frac{1}{2}x(x-2)^2, & \text{for } 1 \leq x \leq 2, \\ N_4(4-x), & \text{for } 2 \leq x \leq 4. \end{cases}$$

Then $\phi_4(k) = a_{k-1}N_4(1) + a_{k-2}N_4(2) + a_{k-3}N_4(3)$ which gives

$$36 \sum \phi_4^2(k) = \sum (a_{k-1} + 4a_{k-2} + a_{k-3})^2.$$

Also, for $t \in [0,1]$,

$$\phi_4(k+t) = a_k N_4(t) + a_{k-1}N_4(t+1) + a_{k-2}N_4(t+2) + a_{k-3}N_4(t+3),$$

and since $N_4(t+2) = N_4(2-t)$ and $N_4(t+3) = N_4(1-t)$, we get

$$\phi_4^2(k+t) = a_k \frac{t^3}{6} + a_{k-1}\left[\frac{2}{3} - \frac{1}{2}(t+1)(t-1)^2\right]$$
$$+ a_{k-2}\left[\frac{2}{3} - \frac{1}{2}(2-t)t^2\right] + a_{k-3}\frac{(1-t)^3}{6}$$
$$= \frac{1}{6}[(a_k - 3a_{k-1} + 3a_{k-2} - a_{k-3})t^3 + (3a_{k-1} - 6a_{k-2} + 3a_{k-3})t^2$$
$$+ 3(a_{k-1} - a_{k-3})t + (a_{k-1} + 4a_{k-2} + a_{k-3})].$$

Then, for $0 \leq t \leq 1$,

$$36 \sum \phi_4^2(k+t)$$
$$= \sum [a_k - 3a_{k-1} + 3a_{k-2} - a_{k-3}]^2 t^6$$
$$+ \sum 2(a_k - 3a_{k-1} + 3a_{k-2} - a_{k-3})(3a_{k-1} - 6a_{k-2} + 3a_{k-3})t^5$$
$$+ \sum [(3a_{k-1} - 6a_{k-2} + 3a_{k-3})^2$$
$$+ 6(a_k - 3a_{k-1} + 3a_{k-2} - a_{k-3})(a_{k-1} - a_{k-3})]t^4$$
(4.52)
$$+ \sum [2(a_k - 3a_{k-1} + 3a_{k-2} - a_{k-3})(a_{k-1} + 4a_{k-2} + a_{k-3})$$
$$+ 6(3a_{k-1} - 6a_{k-2} + 3a_{k-3})(a_{k-1} - a_{k-3})]t^3$$

$$+ \sum [9(a_{k-1} - a_{k-3})^2$$
$$+ 2(3a_{k-1} - 6a_{k-2} + 3a_{k-3})(a_{k-1} + 4a_{k-2} + a_{k-3})]t^2$$
$$+ \sum 6(a_{k-1} - a_{k-3})(a_{k-1} + 4a_{k-2} + a_{k-3})t + 36 \sum \phi_4^2(k).$$

Set
$$A = \sum [a_k - 3a_{k-1} + 3a_{k-2} - a_{k-3}]^2$$
$$= \sum [(a_k - a_{k-1}) - 2(a_{k-1} - a_{k-2}) + (a_{k-2} - a_{k-3})]^2$$
$$= 6 \sum (a_k - a_{k-1})^2 - 8 \sum (a_k - a_{k-1})(a_{k-1} - a_{k-2})$$
$$+ 2 \sum (a_k - a_{k-1})(a_{k-2} - a_{k-3}).$$

Of course, $A \geq 0$, and also, if we set
$$B = -2 \sum (a_k - a_{k-1})(a_{k-1} - a_{k-2}) + 2 \sum (a_k - a_{k-1})(a_{k-2} - a_{k-3}),$$
then,
$$A - B := C = 6 \sum (a_k - a_{k-1})^2 - 6 \sum (a_k - a_{k-1})(a_{k-1} - a_{k-2})$$
$$= 3 \sum ((a_k - a_{k-1}) - (a_{k-1} - a_{k-2}))^2 \geq 0.$$
(4.53)

Actually, $C$ and $A$ are both strictly positive; if $C = 0$ then the points $(k, a_k)$ are on a straight line which is impossible as the $|a_k|$ decrease exponentially with $|k|$. If $A = 0$, then the points $(k, a_k - a_{k-1})$ are on a straight line, and this would also contradict exponential decay of the $a_k$'s [since, for some $c$ and $m$ we would have $a_k = c + mk + a_{k-1} = \cdots = ck + k(k+1)m/2 + a_0$].

Cumbersome but easy manipulation in the above expression for $\sum \phi_4^2(k + t)$ gives the following:

$$36 \sum \phi_4^2(k + t)$$
$$= At^6 - 3At^5 + (2A + B)t^4 + (A - 2B)t^3 + (-A + B)t^2 + \sum \phi_4^2(k)$$
(4.54)
$$= t^2(1-t)^2[(t^2 - t - 1)A + B] + 36 \sum \phi_4^2(k)$$
$$= t^2(1-t)^2[t(t-1)A - C] + 36 \sum \phi_4^2(k).$$

Since $A > 0$ and $C > 0$, we have $t(t-1)A - C < 0$ on $[0, 1]$ and so it follows that $\sigma_4^2(t) = \sum_k \phi_4^2(t + k)$ attains its absolute maximum on $[0, 1]$, namely $\sum \phi_4^2(k)$, only at the points $t = 0$ and $t = 1$. That is, the variance function in the cubic spline-wavelet case attains its absolute maxima, which are strict local maxima, exactly at the points $k \in \mathbb{Z}$.

4.4.4. *A useful inequality.* The following exponential inequality, which is based on Talagrand's (1996) inequality and relevant empirical process techniques, was used repeatedly throughout the proofs.

PROPOSITION 10. *Let $K:\mathbb{R} \to \mathbb{R}$ be either a convolution kernel that is integrable and of bounded variation, or let $K$ be a wavelet projection kernel and assume either that $\phi$ has compact support and is of bounded variation, or that $\phi$ is a Battle–Lemarié father wavelet for some $r \geq 1$. Suppose $P$ has a bounded density $f$ and let $f_n(y,j)$ be the estimator from (2.1). Given $C, T > 0$, there exist finite positive constants $C_1 = C_1(C,K)$ and $C_2 = C_2(C,T,K)$ such that, if*

$$\frac{n}{2^j j} \geq C \quad and \quad C_1 \sqrt{(\|f\|_\infty \vee 1) \frac{2^j j}{n}} \leq t \leq T,$$

*then, for every $n \in \mathbb{N}$,*

$$(4.55) \quad \begin{aligned} \Pr_f \Big\{ \sup_{y \in \mathbb{R}} |f_n(y,j) - Ef_n(y,j)| \geq t \Big\} \\ \leq C_2 \exp\bigg(-\frac{nt^2}{C_2(\|f\|_\infty \vee 1) 2^j}\bigg). \end{aligned}$$

PROOF. Effectively this inequality was proved in Giné and Guillou [(2002), for convolution kernels], Giné and Nickl [(2009b), for compactly supported wavelets] and Giné and Nickl [(2010), for Battle–Lemarié-wavelets]. In the present form it can be deduced, for instance, from Proposition 1 in Giné and Nickl (2010), using the VC-bounds and variance computations in the aforementioned papers, with $\sigma^2 = 2^j c^2(K)(\|f\|_\infty \vee 1)$ and $\lambda$ equal to a large constant (depending only on $C, K$) times $(\|f\|_\infty \vee 1)^{1/2}$. □

**Acknowledgments.** We would like to thank the referees as well as Marc Hoffmann and Aad van der Vaart for very valuable criticism and remarks. We are further grateful to Dominique Picard and Gérard Kerkyacharian for encouraging parts of this research, and to Jürg Hüsler, Wolodomyr Madych as well as Richard Samworth for helpful comments.

DEPARTMENT OF MATHEMATICS  
UNIVERSITY OF CONNECTICUT  
STORRS, CONNECTICUT 06269-3009  
USA  
E-MAIL: gine@math.uconn.edu

STATISTICAL LABORATORY  
DEPARTMENT OF PURE MATHEMATICS  
　AND MATHEMATICAL STATISTICS  
UNIVERSITY OF CAMBRIDGE  
WILBERFORCE ROAD, CB3 0WB, CAMBRIDGE  
UK  
E-MAIL: r.nickl@statslab.cam.ac.uk